%%%%%%%%%%%%%%%%%%%%%%%%%%%%%%%%%%%%%%%%%%%%%%%%%%%%%%%%%%
%%%%%%%%%%%%%%%%%%%%%%%%%%%%%%%%%%%%%%%%%%%%%%%%%%%%%%%%%%
%%
%%     This is the AMS-LaTeX file:
%%
%%     Sprekels-Wu2
%%     Optimal control of a Cahn--Hilliard--Darcy system with mass source
%%     modeling tumor growth
%%
%%%%%%%%%%%%%%%%%%%%%%%%%%%%%%%%%%%%%%%%%%%%%%%%%%%%%%%%%%

\def\LaTeX{%
  \let\Begin\begin
  \let\End\end
  \let\salta\relax
  \let\finqui\relax
  \let\futuro\relax}

\def\UK{\def\our{our}\let\sz s}
\def\USA{\def\our{or}\let\sz z}

\UK
%\USA

%%%%%%%%%%%%%%%%%%%%%%%%%%%%%%%%%

% scegliere fra \TeX e \LaTeX  e fra  \UK oppure \USA

%\TeX
\LaTeX

%\UK
\USA

%%%%%%%%%%%%%%%%%%%%%%%%%%%%%%%%%
%% page layout
%%%%%%%%%%%%%%%%%%%%%%%%%%%%%%%%%

\salta

\documentclass[twoside,12pt]{article}
\setlength{\textheight}{24cm}
\setlength{\textwidth}{16cm}
\setlength{\oddsidemargin}{2mm}
\setlength{\evensidemargin}{2mm}
\setlength{\topmargin}{-15mm}
\parskip2mm

%%%%%%%%%%%%%%%%%%%%%%%%%%%%%%%%%
%% packages
%%%%%%%%%%%%%%%%%%%%%%%%%%%%%%%%%

%\usepackage{color}
\usepackage[usenames,dvipsnames]{color}
\usepackage{amsmath}
\usepackage{amsthm}
\usepackage{amssymb}
\usepackage[mathcal]{euscript}
\usepackage{cite}
%\usepackage[notref,notcite]{showkeys}
%\usepackage{showkeys}
%
%		COLORS FOR CORRECTIONS
%
% do the same, please (i.e., don't use the standard {\color{red} text} or similar):
% just choose the color you prefer in \def\yourname

% example of use:  \hao{I want this to become red}

\def\juerg #1{{\color{red}#1}}

\bibliographystyle{plain}

%%%%%%%%%%%%%%%%%%%%%%%%%%%%%%%%%
%% environments
%%%%%%%%%%%%%%%%%%%%%%%%%%%%%%%%%

%

%finqui

\def\Beq{\Begin{equation}}
\def\Eeq{\End{equation}}

\def\Bthm{\Begin{theorem}}
\def\Ethm{\End{theorem}}
\def\Blem{\Begin{lemma}}
\def\Elem{\End{lemma}}
\def\Bprop{\Begin{proposition}}
\def\Eprop{\End{proposition}}
\def\Bcor{\Begin{corollary}}
\def\Ecor{\End{corollary}}
\def\Brem{\Begin{remark}\rm}
\def\Erem{\End{remark}}

\def\Bdim{\Begin{proof}}
\def\Edim{\End{proof}}
\def\Bcenter{\Begin{center}}
\def\Ecenter{\End{center}}
\let\non\nonumber

%%%%%%%%%%%%%%%%%%%%%%%%%%%%%%%%%
%% macros
%%%%%%%%%%%%%%%%%%%%%%%%%%%%%%%%%

% macro salvate

% sottosezioni non numerate

\def\step #1 \par{\medskip\noindent{\bf #1.}\quad}

% bold, cal e mathop

\def\multibold #1{\def\arg{#1}%
  \ifx\arg\pto \let\next\relax
  \else
  \def\next{\expandafter
    \def\csname #1#1#1\endcsname{{\bf #1}}%
    \multibold}%
  \fi \next}

\def\pto{.}

\def\multical #1{\def\arg{#1}%
  \ifx\arg\pto \let\next\relax
  \else
  \def\next{\expandafter
    \def\csname cal#1\endcsname{{\cal #1}}%
    \multical}%
  \fi \next}

% operatori

\def\multimathop #1 {\def\arg{#1}%
  \ifx\arg\pto \let\next\relax
  \else
  \def\next{\expandafter
    \def\csname #1\endcsname{\mathop{\rm #1}\nolimits}%
    \multimathop}%
  \fi \next}

\multibold
qwertyuiopasdfghjklzxcvbnmQWERTYUIOPASDFGHJKLZXCVBNM.

\multical
QWERTYUIOPASDFGHJKLZXCVBNM.

\multimathop
diag dist div dom mean meas sign supp .

\def\i0t {\int_0^t}
\def\i0T {\int_0^T}
\def\iQt{\int_{Q_t}}
\def\iO{\int_\Omega}
\def\intQ{\int_Q}

\def\vph{\varphi^h}
\def\muh{\mu^h}
\def\uh{{\bf u}^h}
\def\ph{p^h}
\def\yh{y^h}
\def\zh{z^h}
\def\wh{{\bf w}^h}
\def\rh{r^h}
\def\vps{\varphi^*}
\def\ps{p^*}
\def\Rs{R^*}
\def\bv{{\bf v}}
\def\mus{\mu^*}
\def\bus{{\bf u}^*}
\def\dt{\partial_t}
\def\dn{\partial_{\bf n}}
\def\sfrac #1#2{{\textstyle\frac {#1}{#2}}}

\def\checkmmode #1{\relax\ifmmode\hbox{#1}\else{#1}\fi}
\def\aeO{\checkmmode{a.\,e.\, in~$\Omega$}}
\def\aeQ{\checkmmode{a.\,e.\ in~$Q$}}

\def\aeS{\checkmmode{a.\,e.\ on~$\Sigma$}}
\def\ae0T{\checkmmode{a.\,e.\ in~$(0,T)$}}

\def\aat{\checkmmode{for a.\,e.~$t\in(0,T)$}}

% insiemi numerici

\def\erre{{\mathbb{R}}}

\def\enne{{\mathbb{N}}}

% spazi di funzioni a valori vettoriali su [0,T], [0,t], [0,s], [0,+\infty), [\delta,T]

% Come ricordare: in generale i simboli L H W  C da soli per gli spazi su (0,T)
% gli stessi raddoppiati per (0,+\infty)
% aggiunta di t o s al simbolo per (0,t) e (0,s)
% aggiunta di d al simbolo semplice o doppio per intervalli (\delta,T) e (\delta,+\infty)
% il simbolo C e i suoi derivati mettono le quadre anziche' le tonde

% Esempi   \L2V   \L\infty\Vp   \W{1,1}H   \C0H   \LL2V   \CC0\Vp   \Ld2V  \CCdH

% spazi di funzioni su \Omega, \Gamma, Q e \Sigma

\def\Lx #1{L^{#1}(\Omega)}
\def\Hx #1{H^{#1}(\Omega)}

\def\Ldue{\Lx 2}

\def\Huno{\Hx 1}
\def\Hdue{\Hx 2}

% spazi di funzioni su Q e S

% lettere greche

\let\vp\varphi

\def\Uad{{\cal R}_{\rm ad}}

\let\TeXchi\chi                         % new \chi, exactly on the baseline
\newbox\chibox
\setbox0 \hbox{\mathsurround0pt $\TeXchi$}
\setbox\chibox \hbox{\raise\dp0 \box 0 }
\def\chi{\copy\chibox}

% quadratino di fine dimostrazione

% abbreviazioni specifiche del lavoro

\def\Pi{\hat\pi}

\def\etan{\eta_n}
\def\xin{\xi_n}
\def\bvn{\bv_n}
\def\qn{q_n}

\def\bu{{\bf u}}

\let\hat\widehat

\allowdisplaybreaks

%%%%%%%%%%%%%%%%%%%%%%%%%%%%%%
\Begin{document}
%%%%%%%%%%%%%%%%%%%%%%%%%%%%%%%%%

%%%%%%%%%%%%%%%%%%%%%%%%%%%%%%%%%
%% front page
%%%%%%%%%%%%%%%%%%%%%%%%%%%%%%%%%

%
\title{Optimal Distributed Control of a Cahn--Hilliard--Darcy System with Mass Sources}
\author{}
\date{}
\maketitle
\Bcenter
\vskip-1cm
{\large\sc J\"urgen Sprekels$^{(1)}$}\\
{\normalsize e-mail: {\tt sprekels@wias-berlin.de}}\\[.25cm]
{\large\sc Hao Wu$^{(2)}$}\\
{\normalsize e-mail: {\tt haowufd@fudan.edu.cn}}\\[.25cm]
$^{(1)}$
{\small Department of Mathematics}\\
{\small Humboldt-Universit\"at zu Berlin}\\
{\small Unter den Linden 6, 10099 Berlin, Germany}\\[2mm]
{\small and}\\[2mm]
{\small Weierstrass Institute for Applied Analysis and Stochastics}\\
{\small Mohrenstrasse 39, 10117 Berlin, Germany}\\[.45cm]
$^{(2)}$
{\small School of Mathematical Sciences}\\
{\small Key Laboratory of Mathematics for Nonlinear Sciences (Fudan University), Ministry of Education}\\
{\small Shanghai Key Laboratory for Contemporary Applied Mathematics}\\
{\small Fudan University}\\
{\small Han Dan Road 220, Shanghai 200433, China}\\[.3cm]
\Ecenter

\vspace{2mm}
\Begin{abstract}\noindent
In this paper, we study an optimal control problem for a two-dimensional Cahn--Hilliard--Darcy system
with mass sources that arises in the modeling of tumor growth.
The aim is to monitor the tumor fraction in a finite time interval in such a way
that both the tumor fraction, measured in terms of a tracking type cost functional,
is kept under control and minimal harm is inflicted to the patient by administering the control,
which could either be a drug or nutrition. We first prove that the optimal control
problem admits a solution. Then we show that the control-to-state operator
is Fr\'echet differentiable between suitable Banach spaces and derive the first-order
necessary optimality conditions in terms of the adjoint variables and the usual
variational inequality.



\End{abstract}
\salta
\pagestyle{myheadings}
\newcommand\testopari{\sc Sprekels \ --- \ Wu}
\newcommand\testodispari{\sc Optimal control of a Cahn--Hilliard--Darcy system}
\markboth{\testodispari}{\testopari}
\finqui
%
%%%%%%%%%%%%%%%%%%%%%%%%%%%%%%%%%
%% very beginning
%%%%%%%%%%%%%%%%%%%%%%%%%%%%%%%%%

\section{Introduction}
\label{Intro}
\setcounter{equation}{0}
In recent years, the study of tumor growth has attracted a lot of interest.
Various types of mathematical models have been developed to capture the dynamics of morphological changes of a growing solid tumor under many effects including cell-cell and cell-matrix adhesion, mechanical stress, cell motility, transport of nutrients, etc. (see \cite{BLM,CL2010,FBG2006,Fri2007,HZO12,Lowen10,Lowen08} and the references therein). These mathematical models will be helpful to understand the complex biological and chemical processes that occur in tumor growth. Moreover, they are important for controlling the spread of a cancerous tumor to
the surrounding tissue and will be helpful to find optimal treatment strategies.

In the classical description, interfaces between the tumor and healthy tissues are usually considered as idealized surfaces of zero thickness, which leads to the so-called sharp interface models. The sharp interface models are often difficult to analyze mathematically, and may break down when the interface undergoes topological changes such as self-intersection,
pinch-off, and splitting. Alternatively, the diffuse interface models replace this hypersurface description of
the interface with a thin layer where microscopic mixing of the macroscopically distinct components of matter are allowed \cite{Anderson,Cahn,LT}.  This not only yields systems of equations that are better amenable to further analysis, but topological changes of the interface can also be handled naturally.

In this paper, we consider the following initial-boundary value problem of the Cahn--Hilliard--Darcy system:
\begin{align}
\label{ssphi}
& \dt\vp\,+\,{\rm div}(\vp\,\bu)\,=\,\Delta\mu\,+\,\widetilde{S} \quad\aeQ,\\[1mm]
\label{ssmu}
&\mu\,=\,-\Delta\vp\,+\,f'(\vp), \,\mbox{ with }\,f(\vp)=\sfrac{1}{4}\,\vp^4-\sfrac{1}{2}\,
\vp^2, \quad\aeQ,\\[1mm]
\label{ssu}
&\bu\,=\,-\nabla p\,+\,\mu\,\nabla\vp \quad\aeQ,\\[1mm]
\label{ssdivu}
&{\rm div}(\bu)\,=\,S\quad\aeQ,\\[1mm]
\label{ssbc}
&\dn\vp=\dn\mu=0 \quad\mbox{and }\,\bu\cdot{\bf n}=0, \quad\aeS,\\[1mm]
\label{ssini}
&\vp(0)=\vp_0 \quad\aeO,
\end{align}
Here, $\Omega\subset\erre^2$ is assumed to be a smooth bounded domain with boundary $\Gamma:=\partial\Omega$. We denote the outward unit normal field by ${\bf n}$, and by $\dn$ the outward normal derivative. Let $T>0$ be a prescribed final time. We set
\begin{align*}
&Q_t:=\Omega\times (0,t),\quad \text{and}\quad \Sigma_t:=\Gamma\times (0,t),\quad \text{for every}\ t\in (0,T),\\
&Q:=\Omega\times (0,T],\quad \text{and}\quad \Sigma:=\Omega\times (0,T].
\end{align*}
The Cahn--Hilliard--Darcy system \eqref{ssphi}--\eqref{ssini} (also referred to as Cahn--Hilliard--Hele--Shaw system in the context of a multi-phase fluid mixture confined in a porous medium like a Hele--Shaw cell),  constitutes a diffuse
interface model that arises in the study of morphological evolution in solid tumor growth,
where all of the relevant physical parameters (including the thickness of the transition
layers) are normalized to unity. It can be viewed as the simplest version of the general thermodynamically consistent diffuse interface model for tumor growth that was derived in \cite{CWSL, Lowen10,Lowen08} based on the principle of mass conservation together with the second law of thermodynamics (see also \cite{GLSS,GLNS} for recent developments). In this continuum framework, the scalar order parameter $\,\vp\,$ represents the tumor volume fraction, $\,\bu\,$ is the advective velocity field, and $\mu\,$ is the chemical potential.
The phase function $\vp$ satisfies a mass balance law that is governed by a Cahn--Hilliard type equation with additional source term, while the cell velocity $\bu$ fulfills a generalized Darcy law where, besides the pressure gradient,
there appears also the so-called Korteweg force due to the cell concentration. We remark that the chemical potential $ \,\mu\,$ is the variational
derivative of the adhesion energy functional
\Beq
\label{free}
{\cal E}(\vp):=\iO\left(\frac{1}{2}|\nabla\vp|^2+f(\vp)\right)dx,
\Eeq
in which the function $f$ given in \eqref{ssmu} can be regarded as a smooth
double-well approximation of the physically relevant logarithmic
potential \cite{Cahn}. The source terms $\,\widetilde{S}\,$ and $\,S\,$ stand for the possible inter-component
mass exchange as well as gains due to proliferation of cells and loss due to cell death.
They may take different forms according to specific considerations in tumor modelling \cite{CL2010,GLSS,HZO12,Lowen08}.
For the sake of simplicity, we assume here that these source terms are given scalar functions only depending on time $t$ and space $x$.
Concerning the scalar function $p$, we see that it satisfies for almost every $t\in(0,T)$ the elliptic
boundary value problem
\begin{align}
\label{ssp}
-\Delta p(t)=S(t)-{\rm div}(\mu(t)\nabla\vp(t))\quad\aeO,\quad\dn p(t)=0\quad\mbox{a.e. on $\,\Gamma$}.
\end{align}
Clearly, $p$ is only determined up to a constant due to the homogeneous Neumann boundary condition. Thus, we make $p$ unique by always requiring that
\Beq
\label{mean}
{\rm mean}\,p(t)\,:=\,\frac{1}{|\Omega|}\int_
\Omega p(t)\,dx\,=\,0 \quad\aat,
\Eeq
where $\,|\Omega|\,$ denotes the two-dimensional area of $\,\Omega$.

Besides extensive numerical studies on the Cahn--Hilliard--Darcy system and its variants (see \cite{CWSL,FW2012,GLNS,Wise2010,Wise2011,Lowen08} and the references therein),
rigorous mathematical analysis has been carried out in the recent literature. For the simplest case with a regular potential function $f$ and vanishing source terms such that $\widetilde{S}=S=0$, well-posedness and long-time behavior of the Cahn--Hilliard--Darcy system have been obtained in \cite{LTZ,WW2012,WZ2013} (see also \cite{GioGrWu} for the recent contribution concerning the physically relevant logarithmic potential). Later in \cite{JWZ}, the authors analyzed the system \eqref{ssphi}--\eqref{ssini} with a regular potential function $f$ and nonzero but equal sources $\widetilde{S}=S$.
They proved the existence of global weak solutions and local strong solutions in both two and three dimensions. Moreover, when the spatial dimension is two,
they were able to demonstrate the existence and unique\-ness of global strong solutions as well as the convergence of any global solution to a single equilibrium as $t\to +\infty$.
Quite recently, progresses have been made for more general tumor growth models in \cite{DFRSS,FLRS,GL2016,GL2018}, where the existence of global weak solutions was established under suitable choices of boundary conditions. Besides, we refer to \cite{CGH15,CGRS,FGR} for the analysis of a related model proposed in \cite{HZO12}, where velocities are set to zero and the state variables are reduced to the tumor cell fraction and the nutrient-rich extracellular water fraction. We also mention that there are recent works on the nonlocal version of the Cahn--Hilliard--Darcy system \cite{DG,DGG} and the related Cahn--Hilliard--Brinkman system \cite{BCG,CG}.

Based on the well-posedness result obtained in \cite{JWZ}, we are interested in studying optimal control problems associated with the system \eqref{ssphi}--\eqref{ssini} at least in the two-dimensional case. For this purpose,
we write throughout this paper
\begin{align}
\widetilde{S}=S+R.\label{R}
\end{align}
The quantity $\,R$ will be taken as our control variable, and it represents an external source (say, a drug or a nutrient) that can be supplied to the system to monitor the size of the tumor fraction $\vp$. Then the optimal control problem under investigation reads as follows:

\vspace{2mm}
\noindent {\bf (CP)} \quad Minimize the tracking type cost functional
\begin{align}
\label{cost}
{\cal J}(\vp,R)&:=\frac{\beta_1}2 \iO |\vp(T)-\vp_\Omega|^2\,dx \,+\,\frac{\beta_2}2
\intQ |\vp-\vp_Q|^2\,dx\,dt\,+\,\frac{\beta_3}{2} \intQ |R|^2\,dx\,dt
\end{align}
subject to the Cahn--Hilliard--Darcy system \eqref{ssphi}--\eqref{ssini}
and to the control constraint
\Beq
\label{uad}
R\in \Uad:=\{R\in L^2(Q):\,R_{\rm min}\,\le\,R\,\le\,R_{\rm max}\quad\aeQ\}.
\Eeq
Here, $R_{\rm min},R_{\rm max}\in L^\infty(Q)$, the nonnegative constants
$\,\beta_1,\beta_2,\beta_3$, and the target functions $\,\vp_\Omega\in L^2(\Omega)$, $\vp_Q\in L^2(Q)\,$ are prescribed.
Then it easily follows that the set of admissible controls $\Uad$ is a \textit{nonempty, closed, bounded and convex} subset in $L^2(Q)$.

\vspace{1mm}
Let us now state the main results of this paper:
\begin{itemize}
\item[(1)] We establish the existence of optimal controls for problem {\bf (CP)} (see Theorem \ref{opexe});
\item[(2)] We show that the control-to-state operator $\mathcal{S}$ defined by the Cahn--Hilliard--Darcy system \eqref{ssphi}--\eqref{ssini} is
Fr\'chet differentiable between appropriate Banach spaces (see Theorem \ref{diffS});
\item[(3)] We derive the first-order necessary optimality condition (see Theorem \ref{noc}),
in particular, in terms of a variational inequality involving the adjoint state (see Corollary \ref{noc1}).
\end{itemize}

To the best of our knowledge, the optimal control problem associated with the Cahn--Hilliard--Darcy system has never been considered in the literature.
Concerning the optimal distributed/boundary control problems for Cahn--Hilliard type systems subject to various boundary conditions, we refer to \cite{CFGS1,CGS3,CGRS1,CGS1,GLR,HiWe12,ZL1,ZL2} and references therein. When the fluid interaction is taken into account, we refer to \cite{RS} for the control problem of a nonlocal convective Cahn--Hilliard system and to \cite{CGS4,CGS5} for a local convective Cahn--Hilliard system with dynamic boundary conditions in three dimensions of space. In both cases, the fluid velocity is used as the control. If one further assumes that the fluid velocity is governed by the Navier--Stokes system, in \cite{FRS}, the distributed optimal control problem of a two-dimensional nonlocal Cahn--Hilliard--Navier--Stokes system was analyzed, and the result was recently extended to the case with degenerate mobility and singular potential in \cite{FGS}. We also would like to mention the papers \cite{HKW17,HiWe14,HiWe17}, which deal with the optimal control problem for the time-discretized version of Cahn--Hilliard/Navier--Stokes systems in three dimensions.

 As far as the Cahn--Hilliard--Darcy system is concerned, we recall that the uniqueness of global weak solutions to system \eqref{ssphi}--\eqref{ssini} with a regular potential still remains to be an open issue
 even when the spatial dimension is two \cite{JWZ,WZ2013}. However, in order to deduce a well-defined control-to-state operator,
 it requires the unique solvability of the state system itself. On the other hand, the derivation of the Fr\'echet differentiability
 of the control-to-state operator also requires that the solution to the state system be sufficiently regular.
 As a consequence, in this paper we have to confine ourselves to the spatially two-dimensional case,
 in which the existence of a unique global strong solution (at least on an arbitrary given time interval $[0,T]$) can be established.
 Besides, we notice that our state system \eqref{ssphi}--\eqref{ssini} indeed differs from the one that has been considered in \cite{JWZ},
 in which the source $S$ was assumed to have zero mean value (in order to be consistent with the homogeneous Neumann boundary conditions \eqref{ssbc}).
 In contrast to this, the control function $R$ can no longer reasonably be expected to have a zero mean value,
which in turn entails that in our situation the mean value of the order parameter $\vp$ is not necessarily conserved in time.
 Finally, we remark that based on our recent work \cite{GioGrWu}, it is possible to extend the results in this paper to the more physical relevant case with singular potentials,
 and this will be illustrated in a forthcoming paper.

The remaining part of this paper is organized as follows: in Section 2, we give the general
setting and present some preliminary results concerning the solutions to the state system \eqref{ssphi}--\eqref{ssini}.
In Section 3, it is shown that the control-to-state mapping is Fr\'echet differentiable between suitable Banach spaces.
The final Section 4 then brings the main results for the control problem {\bf (CP)}: existence of optimal controls, and the first-order
necessary optimality conditions in terms of a variational inequality and the adjoint variables.

%%%%%%%%%%%%%%%%%%%%%%%%%%%%%%%%%%%%%%%%%%%%%%%%%%%%%%%%%%%%%%%%%%%%%%%%

\section{Results for the State System}
\setcounter{equation}{0}

In this section, we introduce some notations and present some results on the state system \eqref{ssphi}--\eqref{ssini}.

\subsection{Preliminaries}
As in the introduction, $\Omega\subset\erre^2$ is assumed to be a smooth bounded domain with boundary $\Gamma:=\partial\Omega$, and we write $|\Omega|$ for its Lebesgue measure.
Throughout the paper, we denote for a general real Banach space $X$ by $X'$ its dual
and by $\langle\cdot\,,\,\cdot\rangle_{X',X}$ the dual
pairing between elements of $X'$ and $X$. Moreover, $\|\cdot\|_X$ stands for the
norm of $X$ or any power of it.
Notice that in two dimensions of space we have the dense, continuous
and compact embeddings $\Huno\subset L^p(\Omega)$, for all $p\in [1,+\infty)$, and
$\Hdue\subset C^0(\overline\Omega)$. In particular, it holds that
\begin{align}
\label{embed1}
\|\vp\|_{L^p(\Omega)}\,&\le\,C_{p,\Omega}\,\|\vp\|_{\Huno} \quad  \forall\,\vp\in \Huno, \quad p\in [1,+\infty),
\end{align}
with positive constants $C_{p,\Omega}$ that depend only on
$\Omega$ and $p$.

We make frequent use of the following special \juerg{cases} of the Ladyshenskaya and Agmon inequalities,
which are valid for regular two-dimensional domains:
\begin{align}
\label{Lady}
\|\vp\|_{L^4(\Omega)}\,&\le\,C_\Omega\,\|\vp\|_{L^2(\Omega)}^{1/2}\,\|\vp\|_{\Huno}^{1/2}\,\quad\forall\,\vp\in
\Huno,\\[1mm]
\label{Agmon}
\|\vp\|_{L^\infty(\Omega)}\,&\le\,C_\Omega\,\|\vp\|_{L^2(\Omega)}^{1/2}\,\|\vp\|_{\Hdue}^{1/2}\quad\,\forall\,
\vp\in \Hdue.
\end{align}
Besides, we recall that standard elliptic estimates imply that for any $\vp\in H^4(\Omega)$ with
$\,\dn\vp=\dn(\Delta\vp)=0$ \,on\, $\Gamma$, it holds that
\begin{align}
\label{H2}
\|\vp\|_{\Hdue}\,&\le\,C_\Omega
\left(\|\Delta\vp\|_{L^2(\Omega)}\,+\,\|\vp\|_{L^2(\Omega)}\right),\\[1mm]
\label{H3}
\|\nabla\vp\|_{\Hdue}\,&\le\,C_\Omega
\left(\|\nabla\Delta\vp\|_{L^2(\Omega)}\,+\,\|\vp\|_{H^1(\Omega)}\right),\\[1mm]
\label{H4}
\|\Delta\vp\|_{\Hdue}\,&\le\,C_\Omega
\left(\left\|\Delta^2\vp\right\|_{L^2(\Omega)}\,+\,\|\Delta\vp\|_{L^2(\Omega)}\right).
\end{align}
In all of the above inequalities \eqref{Lady}--\eqref{H4}, the positive constant $C_\Omega$ depends only on $\Omega$.

\subsection{Well-posedness and continuous dependence with respect to the source term}

We shall make the following assumptions on the initial datum and source terms of our system.

\vspace{1.5mm}\noindent
(A1) \quad$\vp_0\in \Hdue$\, and \,$\dn\vp_0=0$ \,a.\,e. on \,$\Gamma$.

\vspace{1.5mm}\noindent
(A2) \quad$S\in L^2(Q)$ and $\int_\Omega S(t) dx =0$ for \,a.\,e. $t\in (0,T)$.

\vspace{1.5mm}\noindent
(A3) \quad There are an open set ${\cal R}\subset L^2(Q)$ such that $\Uad\subset{\cal R}$, and
a constant $\widehat R>0$ such\\
\hspace*{13mm} that $\,\|R\|_{L^2(Q)}\,\le\,\widehat R\,$ for every $R\in{\cal R}$.

\vspace{2mm}
First, we present the following result on existence and uniqueness of global strong solutions to problem \eqref{ssphi}--\eqref{ssini}.
\Bthm\label{exe2D}
Suppose that the assumptions {\rm (A1)--(A3)} are satisfied. Then for any given $T\in (0,+\infty)$ and $\,R\in {\cal R}$, the state system \eqref{ssphi}--\eqref{ssini} with \eqref{ssp}--\eqref{mean}, admits a unique solution quadruple \,$(\varphi,\mu,\bu,p)$\,  such
 that $\,{\rm mean}\,p(t)=0\,$ for
a.\,e. $t\in (0,T)$ and
\begin{align}
\label{regphi}
\vp&\in H^1(0,T;\Ldue)\cap C^0([0,T];\Hdue)\cap L^2(0,T;H^4(\Omega)),\\
\label{regmu}
\mu&\in C^0([0,T];\Ldue)\cap L^2(0,T;\Hdue),\\
\label{regup}
\bu&\in L^2(0,T;{\Huno}^2),\quad p\in L^2(0,T;\Hdue).
\end{align}
Moreover, there is a constant $M_1>0$, which depends only on $\|\varphi_0\|_{H^2(\Omega)}$, $\Omega$, $T$, $\|S\|_{L^2(Q)}$ and $\,\widehat R$, such that
\begin{align}
\label{ssbounds1}
&\|\vp\|_{H^1(0,T;\Ldue)\cap C^0([0,T];\Hdue)\cap L^2(0,T;H^4(\Omega))}\,+\,\|\mu\|
_{C^0([0,T];\Ldue)\cap L^2(0,T;\Hdue)}\non\\[1mm]
&\quad + \|\bu\|_{L^2(0,T;{\Huno}^2)}\,+\,\|p\|_{L^2(0,T;\Hdue)}\,\le\,M_1\,.
\end{align}
\Ethm
\Bdim
We recall that a similar result on the existence, uniqueness and regularity has been proved
in \cite[Theorem 2.2]{JWZ} for the
special case $\,R=0\,$ by using a suitable Faedo--Galerkin procedure (see also \cite[Theorem 1.1]{LTZ} for the case $S=R=0$).
A closer inspection of the proof given there reveals, however, that only minor and straightforward modifications are necessary to get these results also in our situation with the nonzero external source term $\,R\in {\cal R}$. Moreover, the uniform estimate \eqref{ssbounds1} can be obtained following the argument in \cite[Lemma 4.1]{JWZ} (see also \cite[Section 4]{LTZ}).
\Edim
\vspace{1mm}
\Brem
(1) In our current setting, we do not need to assume the zero mean constraint for the source term $R$. However, this leads to the fact that the total mass $\int_\Omega \varphi(t) dx$ is no longer conserved as in \cite{JWZ}. On the other hand, one can easily deduce that
\begin{align}
\left|\int_\Omega\varphi(t) \,dx\right| &=\left| \int_\Omega\varphi_0\, dx+ \int_0^t\int_\Omega R(s) \,dx\,ds\right|\non\\
&\leq \|\varphi_0\|_{L^1(\Omega)}+T^{1/2}|\Omega|^{1/2}\widehat{R},\quad \forall\, t\in [0,T].\nonumber
\end{align}

(2) By virtue of \eqref{regphi}--\eqref{ssbounds1} and the Sobolev embedding theorem in two dimensions,
we may without loss of generality assume that
\Beq
\label{ssbounds2}
\max_{0\le j\le 3}\,\left\|f^{(j)}(\vp)\right\|_{C^0([0,T];L^\infty(\Omega))}\,+\,\|\mu\|_{L^2(0,T;
L^\infty(\Omega))}\,+\,\|p\|_{L^2(0,T;L^\infty(\Omega))}\,\le\,M_1,
\Eeq
by choosing $M_1$ properly larger.
\Erem

The next result is concerned with the stability property of the strong solution (i.e., continuous dependence) with respect to the control parameter $R$ in a suitable topology.
It will play an important role in proving the differentiability property of the control-to-state mapping defined by the state problem \eqref{ssphi}--\eqref{ssini}.
\Blem\label{conti}
Suppose that the assumptions {\rm (A1)--(A3)} are fulfilled. If $R_i\in {\cal R}$, $i=1,2$, are
given, and $(\vp_i,\mu_i,\bu_i,p_i)$, $i=1,2$, are the corresponding unique strong solutions to the
state system \eqref{ssphi}--\eqref{ssini} which satisfy \eqref{regphi}--\eqref{regup},
then there exists a constant $M_2>0$, which depends only on $\|\varphi_0\|_{H^2(\Omega)}$, $\Omega$, $T$, $\|S\|_{L^2(Q)}$ and $\widehat R$, such that the following estimate holds true:
\begin{align}\label{stabu1}
&\|\vp_1-\vp_2\|_{H^1(0,t;\Ldue)\cap C^0([0,t];\Hdue)\cap L^2(0,t;H^4(\Omega))}\,+\,\|\mu_1-\mu_2\|_{L^2(0,t;\Hdue)}\non\\
&\qquad +\|\bu_1-\bu_2\|
_{L^2(0,t;L^4(\Omega))}+\,\|p_1-p_2\|_{L^2(0,t;W^{2,4/3}(\Omega))}\non\\
&\quad \le\,M_2\,\|R_1-R_2\|_{L^2(0,t;\Ldue)},\qquad \forall\, t\in (0,T].
\end{align}
\Elem
\Bdim
Let $$R=R_1-R_2,\quad \vp=\vp_1-\vp_2,\quad \mu=\mu_1-\mu_2,\quad \bu=\bu_1-\bu_2,\quad p=p_1-p_2.$$
It then readily follows that
\begin{align}
\label{diff1}
\dt\vp-\Delta\mu\,&=\,R-S\,\vp-\bu\cdot\nabla\vp_1-\bu_2\cdot\nabla\vp\quad\aeQ,\\[1mm]
\label{diff2}
\mu\,&=\,-\Delta\vp+f'(\vp_1)-f'(\vp_2)\quad\aeQ,\\[1mm]
\label{diff3}
\bu\,&=\,-\nabla p+\mu\,\nabla\vp_1+\mu_2\,\nabla\vp\quad\aeQ,\\[1mm]
\label{diff4}
{\rm div}(\bu)\,&=\,0\quad\aeQ,\\[1mm]
\label{diff5}
-\Delta p\,&=\,-\,{\rm div}(\mu\,\nabla\vp_1)-{\rm div}(\mu_2\,\nabla\vp)\quad\aeQ,\\[1mm]
\label{diff6}
\dn\vp\,&=\,\dn\mu\,=\,\dn p\,=\,\bu\cdot{\bf n}\,=\,0\quad\aeS,\\[1mm]
\label{diff7}
\vp(0)\,&=\,0\quad\aeO\,.
\end{align}
Moreover, we have $\,{\rm mean}\, p(t)=0\,$ for almost every $t\in (0,T)$. In what follows, we shall simply perform the estimates in a formal manner, since the arguments can be made rigorous within the approximation scheme devised in \cite{JWZ}.

Let
$t\in (0,T]$ be fixed, but arbitrary. In addition, the letters \,$C$\,
and $\,C_i$, $i\in\enne$, denote positive
constants that only depend on the data of the state system \eqref{ssphi}--\eqref{ssini}
and  on $\widehat R$. The actual meaning
of \,$C$\, may change between or even within lines.
We also notice that \eqref{diff2}, \eqref{diff6} and the uniform estimate \eqref{ssbounds2} imply that
\begin{align}
\label{diff14}
&|{\rm mean}\,\mu(t)|\,=\,|\Omega|^{-1}\iO|f'(\vp_1(t))-f'(\vp_2(t))|\,dx\,\le\,C\,\|\vp(t)\|_{L^2(\Omega)}\,.
\end{align}
Therefore, it follows from Poincar\'e's inequality that
\Beq\label{diff15}
\int_0^t\|\mu(s)\|_{\Huno}^2\,ds\,\le\,C_1\int_0^t(\|\nabla\mu(s)\|^2_{\Ldue}\,+\,\|\vp(s)\|_{\Ldue}^2)\,ds\,.
\Eeq

\vspace{2mm}\noindent
{\bf First estimate:}

\vspace{1mm}\noindent
We multiply \eqref{diff1} by $\vp$, account for \eqref{diff2}, and integrate over $Q_t$. After integration by parts, we obtain that
\begin{align}
\label{diff9}
& \frac 12\,\|\vp(t)\|_{\Ldue}^2\,+\iQt|\Delta\vp|^2\,dxds\non\\
&\quad =\,\iQt\Delta\vp\,(f'(\vp_1)-f'(\vp_2))\,dxds \,+\iQt R\vp\, dxds\,-\iQt
\vp\,(\vp \,S\,+\,\bu_2\cdot\nabla\vp)\,dxds\non\\
&\qquad -\iQt\vp\,(\bu\cdot\nabla\vp_1)\, dxds\non\\
&\quad =\,\sum_{j=1}^4I_j\,,
\end{align}
with obvious notation. We have, for any $\gamma\in (0,1)$ (to be chosen later),
\begin{align}
\label{diff10}
|I_1|+|I_2|\,&\le\,\gamma\iQt|\Delta\vp|^2\, dxds \,+\,C\iQt|R|^2\, dxds \,+\,C(1+\gamma^{-1})\iQt|\vp|^2\, dxds\,,\\
\label{diff11}
|I_3|\,&\le\,\int_0^t\|\vp(s)\|_{L^4(\Omega)}\left(\|\bu_2(s)\|_{L^4(\Omega)}\,\|\nabla\vp(s)\|_{\Ldue}
\,+\,\|\vp(s)\|_{L^4(\Omega)}\,\|S(s)\|_{\Ldue}\right)ds\non\\
&\le\,C\int_0^t\left(\|\bu_2(s)\|_{L^4(\Omega)}\,+\,\|S(s)\|_{\Ldue}\right)\|\vp(s)\|_{\Huno}^2\,ds\,,\\
\label{diff12}
|I_4|\,&\le\int_0^t\|\vp(s)\|_{L^4(\Omega)}\,\|\bu(s)\|_{\Ldue}\,\|\nabla\vp_1(s)\|_{L^4(\Omega)}\,ds\non\\[1mm]
&\le\,\frac 14\,\int_0^t \|\bu(s)\|_{\Ldue}^2\,ds +\,C\int_0^t\|\vp(s)\|_{\Huno}^2\, ds,
\end{align}
where in the last estimate we employed the uniform estimate \eqref{ssbounds1}.
Combining \eqref{diff9}--\eqref{diff12}, we thus have obtained the estimate
\begin{align}
\label{diff13}
&\frac 12\,\|\vp(t)\|_{\Ldue}^2\,+\left(1-\gamma\right)
\iQt|\Delta\vp|^2\,dxds \non\\
&\quad \le\, \frac14 \,\int_0^t \|\bu(s)\|_{\Ldue}^2\,ds\,+\,C\iQt|R|^2\, dxds\non\\
&\qquad +\,C\,(1+\gamma^{-1})\int_0^t\left(1+\,\|\bu_2(s)\|_{L^4(\Omega)}\,+\,\|S(s)\|_{\Ldue}\right)\|\vp(s)\|_{\Huno}^2\,ds\,.
\end{align}

\vspace{2mm}\noindent
{\bf Second estimate:}

\vspace{1mm}\noindent
We multiply \eqref{diff1} by $\,\mu$, take the scalar product of \eqref{diff3} with $\bu$, add
the two resulting identities and integrate over $Q_t$. After integration by parts, and in view of \eqref{diff2}, we obtain the identity
\begin{align}
\label{diff16}
&\frac 12\,\|\nabla\vp(t)\|_{\Ldue}^2\,+\iQt\!|\nabla\mu|^2\,dxds\,+\iQt\!|\bu|^2\,dxds\non\\
&\quad =\,\iQt\! R\,\mu \,dxds\,+\iQt\!\mu_2\,(\bu\cdot\nabla\vp)
\,dxds -\iQt\!\mu(S\,\vp\,+\,\bu_2\cdot\nabla\vp)\,dxds \non\\
&\qquad -\iQt(f'(\vp_1)-f'(\vp_2))\,\dt\vp\,dxds\non\\
&\quad =\,\sum _{j=1}^4 J_j\,,
\end{align}
with obvious notation. Clearly, by virtue of \eqref{diff15} and Young's inequality, we get
\begin{align}
|J_1|\,&\le\,\gamma\iQt|\mu|^2\,dxds\,+\,\frac C\gamma\iQt|R|^2\,dxds\non\\
&\le\,\gamma\,C_1\iQt|\nabla\mu|^2\,dxds\,+\,
C\,(\gamma+\gamma^{-1})\iQt(|R|^2\,+\,|\vp|^2)\,dxds.
\label{diff17}
\end{align}
Moreover, by using H\"older's inequality, \eqref{Lady} and \eqref{ssbounds1}, we have that
\begin{align}
\label{diff18}
|J_2|\,&\le\int_0^t\!\|\mu_2(s)\|_{L^4(\Omega)}\,\|\bu(s)\|_{\Ldue}\,\|\nabla\vp(s)\|_{L^4(\Omega)}\,ds\non\\
&\le\,\frac 14\iQt\!|\bu|^2\,dxds\,+\,C\int_0^t\!\|\mu_2(s)\|_{L^2(\Omega)} \|\mu_2(s)\|_{H^1(\Omega)} \|\nabla\vp(s)\|_{\Ldue}\,\|\nabla\vp(s)\|_{\Huno}\,ds \non\\
&\le\,\frac 14\iQt|\bu|^2\,dxds\,+\,\gamma\iQt|\Delta\vp|^2\,dxds\non\\
&\quad\ +\, C(1+\gamma^{-1})\int_0^t(1+\,\|\mu_2(s)\|_{H^1(\Omega)}^2)\|\vp(s)\|_{\Huno}^2\,ds\,.
\end{align}
We also, owing to \eqref{diff15}, have that
\begin{align}
\label{diff19}
|J_3|\,&\le\int_0^t\|\mu(s)\|_{L^4(\Omega)}\left(\|S(s)\|_{\Ldue}\,\|\vp(s)\|_{L^4(\Omega)}\,
+\,\|\bu_2(s)\|_{L^4(\Omega)}\,\|\nabla\vp(s)\|_{\Ldue}\right)ds
\non\\
&\le\int_0^t\|\mu(s)\|_{\Huno}\,(\|S(s)\|_{\Ldue}\,+\,\|\bu_2(s)\|_{L^4(\Omega)})\,\|\vp(s)\|_{\Huno}\,ds\non\\
&\le\,\gamma\iQt|\nabla\mu|^2\,dxds\non\\
&\quad +\, C(\gamma +\gamma^{-1})\int_0^t(1+\|S(s)\|_{\Ldue}^2\,+\,\|\bu_2(s)\|_{L^4(\Omega)}^2)\,\|\vp(s)\|^2_{\Huno}\,ds\,.
\end{align}
For the estimation of $J_4$, we have to substitute for $\dt\vp$ using equation \eqref{diff1}. After integration by parts, we obtain that
\begin{align}
\label{diff20}
J_4\,=&\,\iQt\nabla\mu\cdot[\nabla(f'(\vp_1)-f'(\vp_2))]\,dxds\,-\iQt R\,(f'(\vp_1)-f'(\vp_2))\, dxds\non\\
&\,+\iQt(S\,\vp+\bu_2\cdot\nabla\vp)\,(f'(\vp_1)-f'(\vp_2))\,dxds\non\\
&\,+\iQt(\bu\cdot\nabla\vp_1)(f'(\vp_1)-f'(\vp_2))\, dxds\non\\
=&\,\sum_{j=1}^4 J_{4}^{(j)}\,,
\end{align}
with obvious notation. From \eqref{ssbounds2}, we see that, almost everywhere in $Q_t$, it holds
\begin{align*}
|\nabla(f'(\vp_1)-f'(\vp_2))|\,&\le\,|f''(\vp_1)-f''(\vp_2)|\,|\nabla\vp_1|\,+\,|f''(\vp_2)|\,|\nabla\vp|\\
&\le\,C\,(|\nabla\vp_1|\,|\vp|+|\nabla\vp|)\,.
\end{align*}
Then it follows from \eqref{ssbounds1} that
\begin{align}
\label{diff21}
&\int_0^t\|f'(\vp_1(s))-f'(\vp_2(s))\|_{\Huno}^2\,ds\non\\
&\quad \le\,C\int_0^t\left(\|\nabla\vp_1(s)\|_{L^4(\Omega)}^2\,\|\vp(s)\|_{L^4(\Omega)}^2
\,+\,\|\vp(s)\|_{\Huno}^2\right)ds\non\\
&\quad \le\,C\int_0^t\|\vp(s)\|_{\Huno}^2\,ds\,.
\end{align}
We thus have, by H\"older's inequality,
\begin{align}
\label{diff22}
|J_{4}^{(1)}|\,&\le\,\gamma\iQt|\nabla\mu|^2\,dxds +\,\frac C\gamma\int_0^t\|f'(\vp_1(s))-f'(\vp_2(s))\|^2_{\Huno}\,ds\non\\
&\le\,\gamma\iQt|\nabla\mu|^2\,dxds +\,\frac C\gamma\int_0^t\|\vp(s)\|_{\Huno}^2\,ds\,.
\end{align}
Moreover, we have
\begin{align}
\label{diff23}
|J_{4}^{(2)}|&\leq\, \int_0^t \|R(s)\|_{\Ldue}\|f'(\vp_1)-f'(\vp_2)\|_{\Ldue}\,ds\non\\
&\le\,C\iQt(|R|^2\,+\,|\vp|^2)\,dxds
\end{align}
as well as
\begin{align}
\label{diff24}
|J_{4}^{(3)}|\,&\le\,\int_0^t\|S(s)\|_{\Ldue}\,\|\vp(s)\|_{L^4(\Omega)}\,
\|f'(\vp_1(s))-f'(\vp_2(s))\|_{L^4(\Omega)}\,ds\non\\
&\quad +\,\int_0^t\|\bu_2(s)\|_{L^4(\Omega)}\,\|\nabla\vp(s)\|_{\Ldue}\,
\|f'(\vp_1(s))-f'(\vp_2(s))\|_{L^4(\Omega)}\,ds\non\\
&\le\,C\int_0^t(1+\|S(s)\|_{\Ldue}^2\,+\,\|\bu_2(s)\|_{L^4(\Omega)}^2)\,\|\vp(s)\|^2_{\Huno}\,ds\,.
\end{align}
Finally, by \eqref{ssbounds1} and \eqref{diff21}, we obtain
\begin{align}
\label{diff25}
|J_{4}^{(4)}|\,&\le\int_0^t\|\bu(s)\|_{\Ldue}\,\|\nabla\vp_1(s)\|_{L^4(\Omega)}\,\|f'(\vp_1(s))-f'(\vp_2(s))\|_{L^4(\Omega)}\,ds\non\\
&\le\,\frac 14\iQt|\bu|^2\,dxds\,+\,C\int_0^t\|\vp(s)\|^2_{\Huno}\,ds\,.
\end{align}
Combining the estimates \eqref{diff13}--\eqref{diff25}, we find that
\begin{align*}
&\frac 12 \,\|\vp(t)\|^2_{\Huno}\,+\,
\frac 14\iQt|\bu|^2\,dxds+\,(1-2\gamma)\iQt|\Delta\vp|^2\,dxds\non\\
&\quad +\,\big(1-(2+C_1)\gamma\big)\iQt|\nabla\mu|^2\,dxds \non\\
&\le\,C(1+\gamma+\gamma^{-1})\iQt|R|^2\,dxds\,+\,C\left(1+\gamma+\gamma^{-1}\right)\int_0^t\Psi_1(s)\,\|\vp(s)\|^2_{\Huno}\,ds\,,
\end{align*}
where, thanks to (A2) and \eqref{ssbounds1}, the function
$$\Psi_1(s)=1+\|\mu_2(s)\|_{H^1(\Omega)}^2+\|S(s)\|_{\Ldue}^2+\|\bu_2(s)\|_{L^4(\Omega)}^2$$
is known to belong to $L^1(0,T)$.
Hence, choosing $\,\gamma\in (0,1/(2+C_1))$, using standard elliptic estimates,
and recalling \eqref{diff15}, we conclude from Gronwall's lemma that
\begin{align}
\label{diff26}
&\|\vp\|_{C^0([0,t];\Huno)\cap L^2(0,t;\Hdue)}^2\,+\,\|\mu\|_{L^2(0,t;\Huno)}^2\,+\,\|\bu\|_{L^2(0,t;\Ldue)}^2\non\\
&\quad \le\,C \iQt|R|^2\,dxds\,.
\end{align}
From \eqref{diff2}, \eqref{diff21}, \eqref{diff26}, we also have
\begin{align}
\label{diff26a}
&\|\vp\|_{L^2(0,t;H^3(\Omega))}^2\le\,C \iQt|R|^2\,dxds\,.
\end{align}

\vspace{2mm}\noindent
{\bf Third estimate:}

\vspace{1mm}\noindent
First, we infer from \eqref{ssbounds1}, \eqref{diff26}, \eqref{diff26a} and H\"older's inequality  that
\begin{align}
\label{stime1}
&\int_0^t\left(\|\nabla\mu(s)\cdot\nabla\vp_1(s)\|^2_{L^{4/3}(\Omega)}\,+\,\|\mu(s)\,\Delta\vp_1(s)\|_{L^{4/3}(\Omega)}^2\right)ds\non\\
&\quad\  \le\int_0^t\left(\|\nabla\mu(s)\|_{\Ldue}^2\,\|\nabla\vp_1(s)\|_{L^4(\Omega)}^2\,+\,\|\mu(s)\|_{L^4(\Omega)}^2\,\|\Delta\vp_1(s)\|_{\Ldue}^2\right)ds
\non\\
&\quad\ \le\,C\int_0^t\|\mu(s)\|_{\Huno}^2\,\|\vp_1(s)\|^2_{\Hdue}\,ds\non\\
&\quad \ \le\,C\iQt|R|^2\,dxds
\end{align}
and
\begin{align}
\label{stime2}
&\int_0^t\left(\|\nabla\mu_2(s)\cdot\nabla\vp(s)\|^2_{L^{4/3}(\Omega)}\,
+\,\|\mu_2(s)\,\Delta\vp(s)\|^2_{L^{4/3}(\Omega)}\right)ds\non\\
&\quad\ \le\int_0^t\left(\|\nabla\mu_2(s)\|_{\Ldue}^2\,\|\nabla\vp(s)\|_{L^4(\Omega)}^2\,+\,\|\mu_2(s)\|_{L^4(\Omega)}^2\,\|\Delta\vp(s)\|_{\Ldue}^2\right)ds\non\\
&\quad\ \le\,C\int_0^t\|\mu_2(s)\|_{\Huno}^2\,\|\vp(s)\|_{\Hdue}^2\,ds\non\\
&\quad\  \leq\, C\int_0^t\|\mu_2(s)\|_{\Ldue}^2\,\|\vp(s)\|_{\Hdue}^2\,ds\non\\
&\qquad \ \ +C \int_0^t \|\mu_2(s)\|_{\Ldue} \|\Delta \mu_2(s)\|_{\Ldue}\,(\|\nabla \Delta \vp(s)\|_{L^2(\Omega)}\|\nabla \vp(s)\|_{\Ldue}+\|\vp(s)\|_{\Ldue}^2)\,ds\non\\[1mm]
&\quad \ \leq C\, \|\mu_2\|_{C^0([0,t];\Ldue)}^2 \, \|\vp\|_{L^2(0,t;\Hdue)}^2 \non\\[1mm]
&\qquad \ \ +\, C\,t^{1/2}\, \|\mu_2\|_{C^0([0,t];\Ldue)} \,\|\Delta \mu_2\|_{L^2(0,t;\Ldue)} \,\|\vp\|_{C^0([0,t];\Ldue)}^2\non\\[1mm]
&\qquad \ \ +\, C \,\|\mu_2\|_{C^0([0,t];\Ldue)}\,\|\Delta \mu_2\|_{L^2(0,t;\Ldue)}
\,\|\nabla \vp\|_{C^0([0,t];\Ldue)} \, \|\nabla \Delta \vp\|_{L^2(0,t;\Ldue)}  \non\\
&\quad \ \le\,C\iQt|R|^2\,dxds\,.
\end{align}
Now notice that
$$\,\,{\rm div}(\mu\nabla\vp_1+\mu_2\nabla\vp)=\nabla\mu\cdot\nabla\vp_1+\mu\,\Delta\vp_1
+\nabla\mu_2\cdot\nabla\vp+\mu_2\,\Delta\vp.$$
Then the continuity of the embedding $W^{1,4/3}(\Omega)\subset L^4(\Omega)$ in two dimensions of space and standard
estimates for elliptic boundary value problems, applied to equation \eqref{diff5}, yield the chain of  estimates
\begin{align}
\label{stime3}
\|\nabla p\|_{L^2(0,t;L^4(\Omega))}^2\,
&\le\,C\,\|\nabla p\|_{L^2(0,t;W^{1,4/3}(\Omega))}^2\,\le\,C\,
\|p\|_{L^2(0,t;W^{2,4/3}(\Omega))}^2\non\\
& \leq \, C\, \|{\rm div}(\mu\nabla\vp_1+\mu_2\nabla\vp)\|_{L^2(0,t;L^{4/3}(\Omega))}\non\\
& \le\, C\iQt|R|^2\,dxds.
\end{align}
On the other hand, from \eqref{Lady}, \eqref{Agmon}, \eqref{diff26} and \eqref{diff26a}, we can deduce the following
estimate:
\begin{align}
\label{stime4}
&\int_0^t\left(\|\mu(s)\nabla\vp_1(s)\|_{L^4(\Omega)}^2\,+\,\|\mu_2(s)\nabla\vp(s)\|_{L^4(\Omega)}^2\right)ds\non\\
&\le\int_0^t\left(\|\mu(s)\|_{L^8(\Omega)}^2\,\|\nabla\vp_1(s)\|_{L^8(\Omega)}^2\,+\,\|\mu_2(s)\|_{L^4(\Omega)}^2\,\|\nabla\vp(s)\|_{L^\infty(\Omega)}^2\right)ds\non\\[2mm]
&\le \,C \,\|\vp_1\|_{C^0([0,t];H^2(\Omega))}^2  \, \|\mu\|^2_{L^2(0,t;\Huno)}\non\\[1mm]
&\quad +\, C\,\|\mu_2\|_{C^0([0,t];L^2(\Omega))}\,\|\mu_2\|_{L^2(0,t;\Huno)}  \,\|\nabla \vp\|_{C^0([0,t]; \Ldue)}
\,\|\nabla \vp\|_{L^2(0,t;\Hdue)}\non\\
&\le \,C\iQt|R|^2\,dxds\,.
\end{align}
Therefore, it follows from \eqref{diff3}, \eqref{stime3}, and \eqref{stime4}, that
\begin{align}
\label{stime5}
\|\bu\|_{L^2(0,t;L^4(\Omega))}^2\,\le\,C\,\int_{Q_t} |R|^2\, dx ds\,.
\end{align}

\vspace{2mm}\noindent
{\bf Fourth estimate:}

\vspace{1mm}\noindent
Inserting the equation \eqref{diff2} for $\,\mu\,$ into \eqref{diff1}, and testing the resulting identity by
$\Delta^2\vp$, we obtain after integration by parts that
\begin{align}
&\frac 12\|\Delta\vp(t)\|_{\Ldue}^2\,+\,\iQt\left|\Delta^2\vp\right|^2\,dxds\non\\
&\quad = \iQt \Big[\Delta(f'(\vp_1)-f'(\vp_2))\,+\,R\,-\,S\vp\, - \,\bu\cdot\nabla\vp_1\,-\,\bu_2\cdot\nabla\vp\Big]\,\Delta^2 \vp\, dxds\non\\
&\quad \leq \frac12 \iQt\left|\Delta^2\vp\right|^2\,dxds + C \iQt|\Delta(f'(\vp_1)-f'(\vp_2))|^2\,dxds + C\iQt |S\vp|^2\, dxds\non\\
&\qquad  + C\iQt|\bu\cdot\nabla\vp_1\,+\,\bu_2\cdot\nabla\vp|^2\,dxds\,+ \, C \iQt|R|^2\,dxds.
\end{align}
Using the estimates \eqref{ssbounds1}, \eqref{ssbounds2}, \eqref{diff26} and \eqref{stime5}, and
invoking H\"older's inequality, it follows  that
\begin{align*}
&\iQt\left(|\Delta(f'(\vp_1)-f'(\vp_2))|^2\ +\,|S\vp|^2\right)\, dxds \non\\
&\quad \le\, C\int_0^t\|\vp(s)\|_{\Hdue}^2\, ds\,+ \, \int_0^t \|S(s)\|_{\Ldue}^2 \,
\|\vp(s)\|_{L^\infty(\Omega)}^2\, ds\non\\
&\quad \le\, C\int_0^t \left((1+\|S(s)\|_{\Ldue}^2\right)\|\Delta \vp(s)\|_{L^2(\Omega)}^2\, ds\non\\
&\qquad\,  + C\Big(1+\|S\|_{L^2(Q_t)}^2\Big)\iQt|R|^2\, dxds
\end{align*}
and
\begin{align*}
&\iQt|\bu\cdot\nabla\vp_1\,+\,\bu_2\cdot\nabla\vp|^2\,dxds\nonumber\\
&\quad \leq C \int_0^t \left(\|\bu(s)\|_{L^4(\Omega)}^2 \, \|\nabla\vp_1(s)\|_{L^4(\Omega)}^2\,+\, \|\bu_2(s)\|_{L^4(\Omega)}^2 \, \|\nabla\vp(s)\|_{L^4(\Omega)}^2 \right) \,ds\non\\
&\quad \leq C \int_0^t \|\bu(s)\|_{L^4(\Omega)}^2\, ds\,+\,C\,\|\vp\|_{C^0([0,t]; H^1(\Omega))}^2 \int_0^t \|\bu_2(s)\|_{L^4(\Omega)}^2 \, ds\non\\
&\qquad +\,  C\int_0^t \|\bu_2(s)\|_{L^4(\Omega)}^2 \, \|\Delta \vp(s)\|_{L^2(\Omega)}^2 \,ds\non\\
&\quad \leq C\int_0^t \|\bu_2(s)\|_{L^4(\Omega)}^2\,\|\Delta \vp(s)\|_{L^2(\Omega)}^2 \,ds + C \iQt|R|^2\, dxds.
\end{align*}
As a consequence, we obtain that
\begin{align*}
& \frac 12\|\Delta\vp(t)\|_{\Ldue}^2\,+\,\frac 12\iQt\left|\Delta^2\vp\right|^2\, dxds\non\\
&\quad \le C\int_0^t \Psi_2(s)\|\Delta \vp(s)\|_{L^2(\Omega)}^2\, ds + C\iQt|R|^2\, dxds,
\end{align*}
where the function $\Psi_2(s)=1+\|S(s)\|_{\Ldue}^2+ \|\bu_2(s)\|_{L^4(\Omega)}^2$ belongs to $L^1(0,T)$. Thus, it follows from Gronwall's lemma that
\begin{align}
\|\Delta\vp\|_{C^0([0,t]; \Ldue)}^2+ \|\Delta^2 \vp\|_{L^2(0,t; \Ldue)}^2\leq C\iQt|R|^2\, dxds.\label{diff27}
\end{align}
Therefore, invoking \eqref{H4}, \eqref{diff26}, \eqref{diff26a} and \eqref{diff27},
we have the estimate
\begin{align}
\label{stime6}
\|\vp\|^2_{C^0([0,t];\Hdue)\cap L^2(0,t;H^4(\Omega))}\,\le\, C\iQt|R|^2\, dxds\, ,
\end{align}
which together with \eqref{diff2} further yields that
\Beq
\label{stabu2}
\|\mu\|^2_{L^2(0,t;\Hdue)}\,\le\,C\,\iQt|R|^2\, dxds\,.
\Eeq
 Finally, using comparison in equation \eqref{diff1} to estimate $\,\dt\vp$, we obtain
\begin{align}
\label{stime7}
\|\vp\|^2_{H^1(0,t;\Ldue)}\,\le\, C\iQt|R|^2\, dxds\,.
\end{align}
This concludes the proof of the assertion.
\Edim

%%%%%%%%%%%%%%%%%%%%%%%%%%%%%%%%%%%%%%%%%%%%%%%%%%%%%%%%%%%%%%%%%%%%%%%%

\section{Differentiability of the Control-to-State Operator}

\setcounter{equation}{0}
\noindent
Set
\begin{align}
\mathcal{V}:= H^1(0,T;\Ldue)\cap C^0([0,T];\Hdue)\cap L^2(0,T;H^4(\Omega)).\non
\end{align}
On account of the well-posedness result given by Theorem \ref{exe2D}, the \emph{control-to-state operator} given by
\begin{align}
{\cal S}:\, L^2(Q)\to \mathcal{V},\quad    R\in L^2(Q)\mapsto {\cal S}(R):=\vp\in \mathcal{V},
\label{ctos}
\end{align}
is well defined and locally bounded, where $\varphi$ is the unique global strong solution to the state system \eqref{ssphi}--\eqref{ssini} on the time interval $[0,T]$ corresponding to the given initial datum $\varphi_0\in H^2(\Omega)$ with $\partial_\mathbf{n}\vp_0=0$ on $\Gamma$ and to the control $R\in L^2(Q)$. Moreover, as a direct consequence of the stability result given by Lemma \ref{conti}, we obtain the following continuity property for ${\cal S}$:
\Bprop
Suppose that the assumptions {\rm (A1)--(A3)} are satisfied. Then the control-to-state operator \,${\cal S}:R\mapsto \vp$\, is locally Lipschitz continuous as a mapping from $L^2(Q)$ into the space $\mathcal{V}$.
\Eprop

In the remaining part of this section, we aim to establish the Fr\'{e}chet differentiability of the control-to-state operator ${\cal S}:R\mapsto\vp$ in suitable Banach spaces.

\subsection{The linearized system}
To this end, we assume that $\Rs\in{\cal R}$ is fixed and the assumptions (A1)--(A3) are satisfied. Set
\begin{align}
\label{star1}
&\vps:={\cal S}(\Rs), \quad\mus:=-\Delta\vps+f'(\vps), \quad\bus=-\nabla\ps+\mus\nabla\vps,
\end{align}
where the pressure variable $\,\ps(t)\,$ is  a
solution to the elliptic problem, for almost every $t\in (0,T)$,
\Beq
\label{star2}
-\Delta\ps(t)=S(t)-{\rm div}(\mus(t)\nabla\vps(t)) \quad\mbox{in }\,\Omega, \quad
\dn\ps(t)=0 \quad\mbox{on }\,\Gamma.
\Eeq
Clearly, $\ps(t)$ is only determined up to a constant, and we make it unique by requesting that
\Beq
\label{star3}
{\rm mean}\,\ps(t)=0 \quad\aat.
\Eeq
Then, by Theorem \ref{exe2D}, the estimates \eqref{ssbounds1}, \eqref{ssbounds2} are valid for the associated solution $(\vps, \mus, \bus, \ps)$.

For a given function $h\in L^2(Q)$, we then consider the linearized system
\begin{align}
\label{lsphi}
& \dt\xi\,+\,{\rm div}(\vps\,\bv)\,+\,{\rm div}(\xi\,\bus)\,=\,
\Delta\eta\,+\,h \quad\aeQ,\\[1mm]
\label{lsmu}
&\eta\,=\,-\Delta\xi\,+\,f''(\vps)\xi, \,\mbox{ with }\,f''(\vps)=3{\vps}^2-1,
\quad\aeQ,\\[1mm]
\label{lsu}
&\bv \,=\,-\nabla q\,+\,\eta\,\nabla\vps\,+\,\mus\nabla\xi \quad\aeQ,\\[1mm]
\label{lsdivu}
&{\rm div}(\bv)\,=\,0\quad\aeQ,\\[1mm]
\label{lsbc}
&\dn\xi=\dn\eta=0 \quad\mbox{and }\,\bv\cdot{\bf n}=0, \quad\aeS,\\[1mm]
\label{lsini}
&\xi(0)\,=\, 0 \quad\aeO.
\end{align}
The system \eqref{lsphi}--\eqref{lsini} can be easily derived by linearizing the state system \eqref{ssphi}--\eqref{ssini} at $(\vp^*, \mu^*, \bu^*, p^*)$.
We note that (provided the involved quantities are sufficiently
smooth) the pressure variable $\,q\,$ solves for almost every $t\in (0,T)$ the
elliptic boundary value problem
\Beq
\label{ellq}
\Delta q(t)\,=\,{\rm div}(\eta(t)\nabla\vps(t))\,+\,{\rm div}(\mus(t)\nabla\xi(t))\quad
\mbox{in }\,\Omega,\quad \dn q(t)=0\quad\mbox{on }\,\Gamma.
\Eeq
Again, $q$ is only determined up to a constant, and we make it unique by generally
requesting that ${\rm mean}\,q(t)=0$ for almost every $t\in (0,T)$. We also remark that
\eqref{ssu}, \eqref{lsu} and \eqref{lsdivu} imply the identity
\Beq
\label{divtog2}
{\rm div}(\vps\,\bv)+{\rm div}(\xi\,\bus)\,=\,\bv\cdot\nabla\vps\,+\,\bus\cdot\nabla\xi\,+S\xi\,,
\Eeq
and it follows from \eqref{ssbc}, \eqref{lsmu} and \eqref{lsbc} that $\,\dn\Delta\xi=0\,$ a.e. on
$\Sigma$, at least formally.
\smallskip

We have the following well-posedness result.
\Blem\label{line}
Suppose that the assumptions {\rm (A1)--(A3)} are fulfilled and that
$\,\Rs\in{\cal R}\,$ is given and $(\vps,\mus,\bus,\ps)$ are defined as in
\eqref{star1}--\eqref{star3}. Then, for every $\,h\in L^2(Q)\,$, the linearized system \eqref{lsphi}--\eqref{ellq}
admits a unique solution $(\xi,\eta, \bv,q)$ such that
$\,{\rm mean}\,q(t)=0\,$ for almost every $t\in (0,T)$ and
\begin{align}
\label{regxi}
&\xi\in H^1(0,T;\Ldue)\cap C^0([0,T];\Hdue)\cap L^2(0,T;H^4(\Omega)),\\[1mm]
\label{regeta}
&\eta\in C^0([0,T];\Ldue)\cap L^2(0,T;\Hdue),\\[1mm]
\label{regv}
&\bv\in L^2(0,T;L^4(\Omega)),\\[1mm]
\label{regq}
&q\in L^2(0,T;W^{2,4/3}(\Omega))\,.
\end{align}
Moreover, the linear mapping $\,h\mapsto\xi\,$ is continuous as a mapping
between $L^2(Q)$ and the space
$$\mathcal{W}:=C^0([0,T];\Huno)\cap L^2(0,T;H^2(\Omega)).$$
\Elem

\Bdim
We shall use a Faedo--Galerkin approximation scheme for the existence result. To this end, consider the Neumann eigenvalue problem
$\,-\Delta w=\lambda\,w$ \,in\, $\Omega$ with $\dn w=0\,$ on $\Gamma$. It is well known
that there exists a nondecreasing sequence $\{\lambda_j\}_{j\in\enne}$ of eigenvalues, where
$0=\lambda_1<\lambda_2\leq ...$ \,and\, $\lim_{j\to\infty}\lambda_j=+\infty$, and a sequence of associated
eigenfunctions $\{e_j\}_{j\in\enne}$ which form a complete orthonormal basis in $\Ldue$. In
particular, $e_1\equiv |\Omega|^{-1/2}$ is a constant function, and, in view of the
orthogonality, we have that $\,\mbox{mean}\, e_j=0$ for every $j>1$.
For any $n\in\enne$, we introduce the $n$-dimensional space
$\,E_n:={\rm span}\{e_1,\ldots,e_n\}\,$ and consider the problem of finding a function of the form
$$\xin(x,t)=\sum_{j=1}^n g_{nj}(t)e_j(x),\quad\mbox{for }\,\, (x,t)\in \overline Q,$$
for some $g_{nj}\in H^1(0,T)$, $1\leq j\leq n$, which satisfies the approximating system
\begin{align}
\label{gs1}
\dt\xin&\,=\,\Delta\etan\,+\,\Pi_n(h-{\rm div}(\vps\,\bvn)-{\rm div}(\xin\,\bus))\quad\aeQ,\\[1mm]
\label{gs2}
\etan&\,=\,-\Delta\xin\,+\,\Pi_n(f''(\vps)\xin), \,\mbox{ with }\,f''(\vps)=3{\vps}^2-1,
\quad\aeQ,\\[1mm]
\label{gs3}
\bvn &\,=\,-\nabla \qn\,+\,\etan \nabla\vps\,+\,\mus\nabla\xin \quad\aeQ,\\[1mm]
\label{gs4}
\xin(0)&\,=\, 0 \quad\aeO.
\end{align}
Here, $\Pi_n$ denotes the $\Ldue$-orthogonal projection onto the space $E_n$, and
$\qn(t)$ solves for almost every $t\in (0,T)$ the elliptic problem
\begin{align}
\label{gs5}
-\Delta \qn(t)\, &=\,-\,{\rm div}(\etan(t)\nabla\vps(t))\,-\,
{\rm div}(\mus(t)\nabla\xin(t))\quad\aeQ,\non\\[1mm]
\dn\qn(t)\,&=\, 0\quad\mbox{a.\,e. on }\,\Sigma,
\end{align}
which entails that
\Beq
\label{gs6}
{\rm div}(\bv_n)=0 \quad\aeQ.
\Eeq
Obviously, $q_n(t)$ is only determined up to a constant. For the sake of simplicity,
and without bearing on the analysis, we require that $\,{\mean}\,q_n(t)=0\,$.

Taking the inner product of \eqref{gs1} and \eqref{gs2} with $\,e_j$, $j=1,\ldots,n$, we
find that the approximating problem is equivalent to an initial value problem for an
explicit system of $\,n\,$ ordinary differential equations for the unknown functions
$\,g_{nj}$, $j=1,\ldots,n$, in which the occurring nonlinearities are locally Lipschitz in the
unknowns and all of the coefficient functions belong to $L^2(0,T)$.
Hence, by virtue of standard results for initial value problems for ordinary differential equations, there
is some $T_n\in (0,T]$ such that the discrete problem \eqref{gs1}--\eqref{gs6} has a unique solution
$(g_{n1},\ldots,g_{nn})\in (H^1(0,T_n))^n$ which specifies the unique local solution of problem \eqref{gs1}--\eqref{gs6} on $\Omega\times(0,T_n)$.
Notice that we have the regularity properties
\begin{align}
\label{reglin}
\xin&\in H^1(0,T_n;\Ldue)\cap L^2(0,T_n;H^4(\Omega)),\non\\[1mm]
\etan&\in H^1(0,T_n;\Ldue)\cap L^2(0,T_n;\Hdue),\non\\[1mm]
\bvn&\in L^2(0,T_n;\Huno), \quad \qn\in L^2(0,T_n;\Hdue),
\end{align}
as well as the boundary conditions
\Beq
\label{bclin}
\dn\xin\,=\,\dn(\Delta\xin)\,=\,\dn\etan\,=\,\dn \qn\,=\,\bv_n\cdot
{\bf n}\,=\,0\quad\mbox{a.\,e. on }\,\Gamma\times (0,T_n).
\Eeq

In the following, we derive some a priori estimates, where the letters \,$C$\,
and $\,C_i$, $i\in\enne$, denote positive constants that may depend on the data of the state system \eqref{ssphi}--\eqref{ssini}, \eqref{ssp}--\eqref{mean}, but neither on $t\in (0,T_n)$ nor on $n\in\enne$. The actual meaning of \,$C$\, may change between or even within lines.
Repeatedly, we shall use the general bounds \eqref{ssbounds1} and \eqref{ssbounds2} for $(\vp^*, \mu^*, \bu^*, p^*)$ without further reference.

Let $t\in (0,T_n)$ be arbitrary, but fixed.
To begin with, we observe that from  \eqref{gs2} and \eqref{bclin} it follows  that
$$|{\rm mean}\, \etan(t)|\,=\,|\Omega|^{-1}\Big|\iO\etan(t)\,dx\Big|\,\le\,C\iO|f''(\vps(t))\,\xin(t)|\,dx
\,\le\,C\,\|\xin(t)\|_{\Ldue}\,,
$$
whence, thanks to Poincar\'e's inequality,
\begin{align}
\label{Poincare}
\int_0^t\|\etan(s)\|^2_{\Huno}\,ds\,\le\,C_1\int_0^t\left(\|\nabla\etan(s)\|_{\Ldue}^2\,
+\,\|\xin(s)\|_{\Ldue}^2\right)\,ds\quad\forall\,n\in\enne.
\end{align}
In what follows, for the sake of a shorter exposition, we
suppress the subscript $n$ during the subsequent calculations, writing it only at the end of the estimations.

\vspace{2mm}\noindent
{\bf First estimate:}

\vspace{1mm}\noindent
We multiply \eqref{gs1} by $\xi_n$ and integrate over $Q_t$.  Using integration by parts, we obtain from \eqref{gs2} the identity
\begin{align}
\label{fe1}
&\frac 12\,\|\xi(t)\|_{\Ldue}^2\,+\iQt\!\! |\Delta\xi|^2\,dxds \non\\
&\quad =\,\iQt\!\! h\,\xi\,dxds\,+\iQt\!\! f''(\vps)\,\xi\,\Delta\xi\,dxds\,-\iQt\!\!{\rm div}(\vps\bv)\,\xi\,dxds\,-\iQt\!\!{\rm div}(\xi\,\bus)\,\xi\,dxds\non\\
&\quad =\,\sum_{j=1}^4 \,I_j\,,
\end{align}
with obvious notation. By virtue of Young's inequality, we have, for every $\gamma\in(0,1)$ (to be chosen later), that
\begin{equation}
\label{fe2}
|I_1|+|I_2|\,\le\,\gamma \iQt\!\!|\Delta\xi|^2\,dxds\,+\,\frac12 \iQt\!\!|h|^2\,dxds\,+\,C\left(1+\gamma^{-1}\right)\iQt\!\!|\xi|^2\,dxds.
\end{equation}
Moreover, from \eqref{divtog2} and H\"older's inequality we conclude that
\begin{align}
\label{fe3}
|I_4|\,&\le\iQt\left(|S||\xi|^2\,+\,|\xi||\bus||\nabla\xi|\right)\,dxds\non\\
&\le \int_0^t \!\|\xi(s)\|_{L^4(\Omega)}\left(\|S(s)\|_{\Ldue}\,\|\xi(s)\|_{L^4(\Omega)}
\,+\,\|\bus(s)\|_{L^4(\Omega)}\,\|\nabla\xi(s)\|_{\Ldue}\right)\,ds
\non\\
&\le\,C\int_0^t\!\!\left(\|S(s)\|_{\Ldue}\,+\,\|\bus(s)\|_{L^4(\Omega)}\right)\|\xi(s)\|^2
_{\Huno}\,ds\,.
\end{align}
Finally, we infer from H\"older's inequality and Young's inequality that
\begin{align}
\label{fe4}
|I_3|\,&\le\,\int_0^t\|\nabla\vps(s)\|_{L^4(\Omega)}\,\|\bv(s)\|_{\Ldue}\,\|\xi(s)\|_{L^4(\Omega)}\,ds\non\\
&\le\,C\int_0^t
\|\bv(s)\|_{\Ldue}\,\|\xi(s)\|_{\Huno}\,ds\non\\
&\le\,\frac 18\int_0^t\|\bv(s)\|_{L^2(\Omega)}^2\,ds\,+\,C\,\int_0^t\|\xi(s)\|_{\Huno}^2\,ds\,.
\end{align}
Combining \eqref{fe1}--\eqref{fe4}, we obtain the following estimate for all $n\in\enne$:
\begin{align}
\label{fe5}
&\,\|\xin(t)\|_{\Ldue}^2\,+\,(2-2\gamma)\iQt|\Delta\xin|^2\,dxds\non\\
&\quad \le\,\iQt|h|^2\, dxds\,+\,\frac 14\int_0^t\|\bv_n(s)\|_{L^2(\Omega)}^2\,ds\non\\
&\qquad +\,C\left(1+\gamma^{-1}\right)\int_0^t\left(1+\|S(s)\|_{\Ldue}+\|\bus(s)\|_{L^4(\Omega)}\right)
\|\xin(s)\|_{\Huno}^2\,ds\,.
\end{align}

\vspace{2mm}\noindent
{\bf Second estimate:}

\vspace{1mm}\noindent
Next, we multiply \eqref{gs1} by $\eta_n$, take the scalar product of \eqref{gs3} with $\bv_n$, invoke \eqref{gs2}, and integrate over $Q_t$. After integration by parts, we obtain the identity
\begin{align}
\label{fe6}
&\frac 12\,\|\nabla\xi(t)\|_{\Ldue}^2\,+\iQt|\nabla\eta|^2\,dxds+\,\iQt|\bv|^2\,dxds\non\\
&\quad =\,\iQt \eta\,(h\,-\,\bus\cdot\nabla\xi\,-\,S\,\xi) \,dxds\,+\,\iQt \mus\,\bv\cdot\nabla\xi\,dxds\,-\,\iQt f''(\vps)\,\xi\,\dt\xi\,dxds\non\\
&\quad :=\,\sum_{j=1}^3 J_j\,,
\end{align}
with obvious notation. Owing to \eqref{Lady}, we infer from H\"older's inequality and Young's inequality that
\begin{align}
\label{fe7}
|J_2|\,&\le\,\int_0^t\|\bv(s)\|_{\Ldue}\,\|\mus(s)\|_{L^4(\Omega)}\,\|\nabla\xi(s)\|_{L^4(\Omega)}\,ds\non\\
&\le\,\frac 14\iQt|\bv|^2\,dxds \,+\,C\int_0^t \|\mus(s)\|_{L^2(\Omega)}\|\mus(s)\|_{H^1(\Omega)}\|\nabla\xi(s)\|_{\Ldue}\,\|\nabla\xi(s)\|_{\Huno}\,ds\non\\
&\le\,\frac 14\iQt|\bv|^2\,dxds\,+\,\gamma\iQt|\Delta\xi|^2\,dxds\non\\
&\quad \ + C\,(1+\gamma^{-1})\int_0^t(1+\|\mus(s)\|_{H^1(\Omega)}^2) \, \|\xi(s)\|^2_{\Huno}\,ds\,.
\end{align}
Moreover, thanks to \eqref{Poincare}, we get
\begin{align}
\label{fe8}
|J_1|\,&\le\gamma\iQt|\eta|^2\,dxds\,+\,\frac C\gamma\iQt|h|^2\,dxds\non\\
&\quad +\int_0^t\|\eta(s)\|_{L^4(\Omega)}\left(\|\bus(s)\|_{L^4(\Omega)}\,\|\nabla\xi(s)\|_{\Ldue}
\,+\,\|S(s)\|_{\Ldue}\,\|\xi(s)\|_{L^4(\Omega)}\right)ds\non\\
& \le\,\gamma\iQt|\eta|^2\,dxds\,+\,\frac C\gamma\iQt|h|^2\,dxds\,+\,C_1\gamma\iQt(|\nabla\eta|^2+|\xi|^2)\,dxds\non\\
&\quad +\,C\gamma^{-1}\int_0^t\left(\|\bus(s)\|_{L^4(\Omega)}^2\,+\,\|S(s)\|_{\Ldue}^2\right)\|\xi(s)\|_{\Huno}^2\,ds\non\\
& \le\,\frac C\gamma\iQt|h|^2\,dxds\,+\,2\,C_1\gamma\iQt |\nabla\eta|^2\,dxds\non\\
&\quad +\,C(1+\gamma+\gamma^{-1})\int_0^t\left(1+\|\bus(s)\|_{L^4(\Omega)}^2\,+\,\|S(s)\|_{\Ldue}^2\right)\|\xi(s)\|_{\Huno}^2\,ds\,.
\end{align}
Finally, we have, using integration by parts in time,
\begin{align}
\label{fe11}
J_3\,&=\,-\frac 12\iO\xi^2(t)f''(\vps(t))\,dx\,+\,\frac 12\iQt\xi^2\,f'''(\vps)\,\dt\vps\,dxds\,:=\,J_4+J_5\,.
\end{align}
It holds that $\,-f''(\vps)=1-3{\vps}^2\le 1$, and thus
\Beq
\label{fe12}
J_4\,\le\,\frac 12\,\|\xi(t)\|_2^2\,.
\Eeq
Besides, it easily follows that
\begin{align}
\label{fe13}
J_5& \le\,C\int_0^t\|\xi(s)\|_{L^4(\Omega)}^2\,\|\dt\vps(s)\|_{\Ldue}\,ds\non\\
&\le\,C\int_0^t\Big(1+\|\dt\vps(s)\|^2_{\Ldue}\Big)\,\|\xi(s)\|_{\Huno}^2\,ds\,.
\end{align}
Combining the estimate \eqref{fe5} and \eqref{fe6}--\eqref{fe13}, we have thus shown that for all $t\in (0,T]$ it holds that
\begin{align}
&\frac 12 \left(\|\xi_n(t)\|_{\Ldue}^2+\|\nabla\xi_n(t)\|_{\Ldue}^2\right)
\,+\,(1-2\,C_1\gamma)\iQt|\nabla\eta_n|^2\,dxds\,\non\\
&\qquad +\,
(2-3\gamma)\iQt|\Delta\xi_n|^2\,dxds\,+\,\frac 12\,\iQt|\bv_n|^2\,dxds\non\\
&\quad \le\,C\left(1+\gamma^{-1}\right)\iQt|h|^2\,dxds\,+\,C\left(1+\gamma+\gamma^{-1}\right)
\int_0^t\Psi_3(s)\,\|\xi_n(s)\|^2_{\Huno}\,ds\,,
\end{align}
where the function  $$\,\,\Psi_3(s):=1+\|S(s)\|_{\Ldue}^2+\|\bus(s)\|^2_{L^4(\Omega)}\,+\,\|\mus(s)\|_{H^1(\Omega)}^2\,+\,\|\dt\vps(s)\|_{\Ldue}^2$$ %
is known to belong to $L^1(0,T)$. Choosing $\gamma\in(0,1)$ appropriately small, and invoking \eqref{Poincare} and standard
elliptic estimates, we then conclude from Gronwall's lemma the estimate
\begin{align}
\label{esti1}
&\|\xin\|^2_{L^\infty(0,t;\Huno)\cap L^2(0,t;\Hdue)}\,+\,\|\etan\|^2_{L^2(0,t;\Huno)}\,+\,\|\bv_n\|^2_{L^2(0,t;\Ldue)}\non\\[1mm]
&\quad \le\,C\,\|h\|^2_{L^2(Q_t)}\quad\mbox{for all $\,t\in (0,T_n)\,$ and $\,n\in\enne.$}
\end{align}
Then we infer from \eqref{gs2}, \eqref{esti1} and the standard
elliptic estimate that
 \begin{align}
\label{esti1a}
&\|\xin\|^2_{L^2(0,t;H^3(\Omega))}\,\le\,C\,\|h\|^2_{L^2(Q_t)}\quad\mbox{for all $\,t\in (0,T_n)\,$ and $\,n\in\enne.$}
\end{align}

\vspace*{2mm}\noindent
{\bf Third estimate:}

\vspace{1mm}\noindent
 Put
\begin{align*}
\psi_n\,:=\,{\rm div}(\etan\nabla\vps)\,+\,{\rm div}(\mus\nabla\xin)\,
=\,\nabla\etan\cdot\nabla\vps+\etan\,\Delta\vps+\nabla\mus\cdot\nabla\xin+\mus\,\Delta\xin.
\end{align*}
Recalling that $\mus\in C^0([0,T];\Ldue)\cap L^2(0,T; \Hdue)$ and $\vps\in C^0([0,T];\Hdue)$, we obtain from \eqref{esti1}, \eqref{esti1a} and H\"older's inequality the chain of estimates (in a similar manner as for \eqref{stime1}, \eqref{stime2})
\begin{align*}
&\int_0^t\|\psi_n(s)\|^2_{L^{4/3}(\Omega)}\,ds\non\\
&\quad \le\,\int_0^t\Big(\|\nabla\etan(s)\|_{\Ldue}^2\,\|\nabla\vps(s)\|_{L^4(\Omega)}^2
\,+\,\|\etan(s)\|_{L^4(\Omega)}^2\,\|\Delta\vps(s)\|_{\Ldue}^2\non\\
&\qquad\quad  +\,\|\nabla\mus(s)\|_{\Ldue}^2\,\|\nabla\xin(s)\|_{L^4(\Omega)}^2\,+\,\|\mus(s)\|_{L^4(\Omega)}^2\,\|\Delta\xin(s)\|_{\Ldue}^2\Big)\,ds\non\\
&\quad \le \,C\int_0^t\left(\|\etan(s)\|_{\Huno}^2\,+\|\xin(s)\|_{\Hdue}^2\right)ds\,+\, C t^{1/2}\|\xin\|_{C^0([0,t]; \Ldue)}^2\non\\
&\qquad +\, C \|\nabla \xin\|_{C^0([0,t]; \Ldue)}\|\nabla \Delta \xin\|_{L^2(0,t;L^2(\Omega))}\non\\
&\quad \le \, C\,\|h\|^2_{L^2(Q_t)}\quad\mbox{for all $\,t\in (0,T_n)\,$ and $\,n\in\enne.$}
\end{align*}
We thus can infer from the standard regularity theory of elliptic boundary value problems
and from the embedding $W^{1,4/3}(\Omega)\subset L^4(\Omega)$, which is valid in two dimensions of space, that
\begin{align}
\label{fe15}
\|\nabla q_n\|^2_{L^2(0,T_n;L^4(\Omega))}
&\le\,C\,\|\nabla q_n\|^2_{L^2(0,T_n;W^{1,4/3}(\Omega))}
\,\le\,C\,\|q_n\|^2_{L^2(0,T_n;W^{2,4/3}(\Omega))}\non\\
&\le\,C\,\|\psi_n\|^2_{L^2(0,T_n;L^{4/3}(\Omega))} \le \, C\,\|h\|^2_{L^2(Q)},\quad\forall\,n\in\enne\,.
\end{align}
Besides, we observe from \eqref{ssbounds1}, \eqref{esti1} and \eqref{esti1a} that, for all $n\in\enne$,
\begin{align}
\label{fe14}
&\int_0^t\|\eta_n(s)\nabla\vps(s)\,+\,\mus(s)\nabla\xi_n(s)\|_{L^4(\Omega)}^2\,ds\non\\
&\quad \le\,C\int_0^t
(\|\eta_n(s)\|_{L^8(\Omega)}^2\,\|\nabla\vps(s)\|_{L^8(\Omega)}^2
\,+\,\|\mus(s)\|_{L^4(\Omega)}^2\,\|\nabla\xi_n(s)\|_{L^\infty(\Omega)}^2)\,ds\non\\
&\quad \le\,C\,\|\vps\|^2_{C^0([0,t];\Hdue)}\int_0^t \|\eta_n(s)\|_{\Huno}^2\,ds\non\\
&\qquad  +\,C\|\mus\|_{C^0([0,t];L^2(\Omega))} \|\nabla\xi_n\|_{C^0([0,t];L^2(\Omega))}
\int_0^t \|\mus(s)\|_{\Huno} \|\nabla \xi_n(s)\|_{\Hdue}\,ds\non\\
&\quad \leq \, C\,\|h\|^2_{L^2(Q_t)}\quad\mbox{for all $\,t\in (0,T_n)\,$ and $\,n\in\enne.$}
\end{align}
Combining \eqref{fe15} and \eqref{fe14}, we finally have shown that
\Beq
\label{f16}
\|\bv_n\|_{L^2(0,T_n;L^4(\Omega))}^2\,\le\,C\,\|h\|^2_{L^2(Q)},\quad\forall\,n\in\enne.
\Eeq

\vspace{2mm}\noindent
{\bf Fourth estimate:}

\vspace{1mm}\noindent
At this point, we are ready to test \eqref{gs1} by $\Delta^2\xin\in E_n$. We obtain, omitting again the
subscript $\,n$, for every $t\in (0,T_n)$ the following identity:
\begin{align}
\label{fe16}
&\frac 12\,\|\Delta\xi(t)\|_{\Ldue}^2\,+\,\iQt\left|\Delta^2\xi\right|^2\,dxds\non\\
&\quad =\,\iQt h\,\Delta^2\xi\,dxds\,+\,
\iQt \Delta(f''(\vps)\,\xi)\,\Delta^2\xi\,dxds\,-\,\iQt(\bv\cdot\nabla\vps)\,\Delta^2\xi\,dxds\non\\
&\qquad -\iQt S\,\xi\,\Delta^2\xi\,dxds-\,\iQt(\bus\cdot\nabla\xi)\,\Delta^2\xi\,dxds\non\\
&\quad :=\,\sum_{j=1}^5K_j\,,
\end{align}
with obvious notation. Using the fact $\vps\in C^0([0,T]; \Hdue)$, Young's inequality and \eqref{esti1}, we have
\begin{align}
\label{fe17}
|K_1|+|K_2|&\le\,\frac 18\iQt\left|\Delta^2\xi\right|^2\,dxds\,+\,C\iQt |h|^2\,dxds\,+\, C\,\int_0^t\|\xi(s)\|_{\Hdue}^2\,ds\non\\
&\le\,\frac 18\iQt\left|\Delta^2\xi\right|^2\,dxds\,+\,C\,\|h\|^2_{L^2(Q_t)}.
\end{align}
Moreover, thanks to \eqref{ssbounds1}, \eqref{f16}, H\"older's inequality and Young's inequality, it holds that
\begin{align}
\label{fe18}
|K_3|&\le\,\int_0^t\|\nabla\vps(s)\|_{L^4(\Omega)}\,\|\bv(s)\|_{L^4(\Omega)}\,\left\|\Delta^2\xi(s)\right\|_{\Ldue}\,ds
\non\\
&\le\,\frac 18\iQt\left|\Delta^2\xi\right|^2\,dxds\,+\,C\,\|h\|^2_{L^2(Q_t)}.
\end{align}
Also, by virtue of \eqref{Lady}, \eqref{regup} and \eqref{esti1}, we obtain
\begin{align}
\label{fe19}
|K_5|\,&\le\int_0^t\|\bus(s)\|_{L^4(\Omega)}\,\|\nabla\xi(s)\|_{L^4(\Omega)}\,\left\|\Delta^2\xi(s)\right\|_{\Ldue}\,ds\non\\
&\le\,\frac 18\iQt\left|\Delta^2\xi\right|^2\,dxds\,+\,C\int_0^t\|\bus(s)\|_{L^4(\Omega)}^2 \,
\|\nabla\xi(s)\|_{\Ldue}\,\|\nabla\xi(s)\|_{\Huno}\,ds
\non\\
&\le\,\frac 18\iQt\left|\Delta^2\xi\right|^2\,dxds\,+\,C\|\xi\|_{C^0([0,t];\Huno)}^2\int_0^t\|\bus(s)\|_{L^4(\Omega)}^2\,ds \non\\
&\quad +\,C\int_0^t\|\bus(s)\|_{L^4(\Omega)}^2\,\|\Delta \xi(s)\|_{\Ldue}^2\,ds\non\\
&\leq \frac 18\iQt\left|\Delta^2\xi\right|^2\,dxds\,+\,C\int_0^t\|\bus(s)\|_{L^4(\Omega)}^2\,\|\Delta \xi(s)\|_{\Ldue}^2\,ds\, +\,C\,\|h\|^2_{L^2(Q_t)}.
\end{align}
Finally, we have
\begin{align}
\label{fe20}
|K_4|\,&\le\,\int_0^t\|S(s)\|_{\Ldue}\,\|\xi(s)\|_{L^\infty(\Omega)}\,\left\|\Delta^2\xi(s)\right\|_{\Ldue}\,ds\non\\
&\le\,\frac 18\iQt\left|\Delta^2\xi\right|^2\,dxds\,+\,C\,\int_0^t\,\|S(s)\|_2^2\,\|\xi(s)\|_{\Ldue}\,\|\xi(s)\|_{\Hdue}\,ds\non\\
&\le\, \frac 18\iQt\left|\Delta^2\xi\right|^2\,dxds\, + \,C\,\int_0^t\,\|S(s)\|_2^2\,\|\Delta \xi(s)\|_{\Ldue}^2\,ds\, +\,C\,\|h\|^2_{L^2(Q_t)}.
\end{align}
Combining \eqref{fe16}--\eqref{fe20}, we obtain that
\begin{align}
\label{fe16a}
&\frac 12\,\|\Delta\xi(t)\|_{\Ldue}^2\,+\,\frac12\, \iQt\left|\Delta^2\xi\right|^2\,dxds\non\\
&\quad \leq C\,\int_0^t\,\Big(\|S(s)\|_2^2+\|\bus(s)\|_{L^4(\Omega)}^2\Big)\,\|\Delta \xi(s)\|_{\Ldue}^2\,ds\,+\,C\,\|h\|^2_{L^2(Q_t)}.
\end{align}
Applying Gronwall's lemma, invoking standard elliptic estimates and \eqref{esti1},  we then conclude  that
\begin{align}
\label{esti4}
\|\xin\|_{L^\infty(0,T_n;\Hdue)\cap L^2(0,T_n;H^4(\Omega))}^2\,\le\,C\,\|h\|^2_{L^2(Q)},\quad\forall\,n\in\enne\,.
\end{align}
From \eqref{gs2} it also easily follows that
\Beq
\label{esti6}
\|\etan\|^2_{L^\infty(0,T_n;\Ldue)\cap L^2(0,T_n;\Hdue)}\,\le\,C\,\|h\|^2_{L^2(Q)},\quad\forall\,n\in\enne.
\Eeq
Moreover, testing \eqref{gs1} by $\dt\xin$, from the above estimates \eqref{esti4}, \eqref{esti6} we can deduce that
\Beq\label{esti5}
\|\xin\|_{H^1(0,T_n;\Ldue)}^2\,\le\,C\,\|h\|^2_{L^2(Q)},\quad\forall   \,n\in\enne.
\Eeq

From the previously shown estimates it immediately follows that $\,T_n=T\,$ for all $\,n\in\enne$.
Furthermore, combining these uniform estimates, we are able to conclude from the well-known compactness results
(cf., in particular,
\cite[Sect.~8,~Cor.~4]{Simon}) that there is some quadruple $(\xi,\eta,\bv,q)$ such that,
at least for some subsequence which is again indexed by $n$ for simplicity,
\begin{align*}
\xin\to\xi\quad&\mbox{weakly-star in \,$H^1(0,T;\Ldue)\cap L^\infty(0,T;\Hdue)\cap L^2(0,T;H^4(\Omega))$} \\
&\mbox{and strongly in }\,C^0(\overline Q),\\[1mm]
\etan\to\eta\quad&\mbox{weakly-star in \,$L^\infty(0,T;\Ldue)\cap L^2(0,T;\Hdue)$,}\\[1mm]
\bv_n\to\bv\quad&\mbox{weakly in \,$L^2(0,T;L^4(\Omega))$,}\\[1mm]
q_n\to q\quad&\mbox{weakly in \,$L^2(0,T;W^{2,4/3}(\Omega))$.}
\end{align*}
In view of these convergence results, it is then a standard matter (which may be left to the reader) to show
that $(\xi,\eta,\bv,q)$ is in fact a solution to the linearized system \eqref{lsphi}--\eqref{ellq} that enjoys the
regularity properties asserted in \eqref{regxi}--\eqref{regq}. Observe, in this connection, that $\xi\in
H^1(0,T;\Ldue)\cap L^2(0,T;H^4(\Omega))\subset C^0([0,T];\Hdue)$, which further entails the continuity of $\eta$, i.e.,
$\eta\in C^0([0,T];\Ldue)$.

\vspace*{2mm}
Next, we show that the solution to the linearized system \eqref{lsphi}--\eqref{ellq} is indeed unique. To this end,
assume that two solution quadruples $(\xi_{i},\eta_{i},\bv_{i},q_{i})$, $i=1,2$, with
the regularity properties \eqref{regxi}--\eqref{regq} are given. Then, for the difference functions
$$\xi= \xi_{1}-\xi_{2},\quad \eta=\eta_{1}-\eta_{2}, \quad \bv=
\bv_{1}-\bv_{2},\quad q=q_{1}-q_{2},
$$
the equations \eqref{lsphi}--\eqref{ellq} are fulfilled with $h=0$. Therefore, the estimations leading to
\eqref{esti1} can be repeated for the continuous system, leading to the conclusion
that \eqref{esti1}  holds true for $(\xi,\eta,\bv)$  with $h=0$, whence it follows that $\xi=\eta=0$ and $\bv=(0,0)$.
As a consequence, we also have $q=0$. Thus, the uniqueness is proved.

It remains to show the continuity of the mapping $h\mapsto \xi$. But this is an immediate consequence
of  \eqref{esti1} if the weak and weak-star sequential lower
semicontinuity properties of the involved norms are taken into account. This concludes the proof
of the assertion.
\Edim

\subsection{Fr\'echet differentiability of $\mathcal{S}$}

We are now in a position to prove the Fr\'echet differentiability of the control-to-state mapping $\mathcal{S}$.
More precisely, we have the following result:
\Bthm\label{diffS}
Suppose that the assumptions {\rm (A1)--(A3)} are fulfilled. Let, for any given
$\,\Rs\in{\cal R}$, the global strong solution $(\vps,\mus,\bus,\ps)$ be defined as in
\eqref{star1}--\eqref{star3}. Then the control-to-state operator ${\cal S}: L^2(Q)\to \mathcal{V}$
defined by \eqref{ctos} is
Fr\'echet differentiable at $\,\Rs\,$ as a mapping from $L^2(Q)$ into the space
${\cal W}:=C^0([0,T];\Huno) \cap L^2(0,T;\Hdue)$. Moreover,
for any $\,h\in L^2(Q)\,$ we have the identity  $$D{\cal S}(\Rs)h=\xi,$$
 where
$(\xi,\eta,\bv,q)$ is the unique solution to the linearized system \eqref{lsphi}--\eqref{ellq} at the point $(\vps, \mus, \bus, \ps)$ that corresponds to the function $h$, subject to the constraint ${\rm mean}\, q(t)=0$ for almost every $t\in (0,T)$.
\Ethm

\Bdim
For any $h\in L^2(Q)$, we denote by $(\xi^h,\eta^h,\bv^h,q^h)$ the unique solution
to the linearized system \eqref{lsphi}--\eqref{ellq} satisfying ${\rm mean}\, q^h(t)=0$ for almost
every $t\in (0,T)$, and we note that by Lemma \ref{line} the linear mapping
$h\mapsto \xi^h$ is continuous between $L^2(Q)$ and the space ${\cal W}$.

Recall that the set
${\cal R}$ is open in $L^2(Q)$ (see (A3)), and thus there is some constant $\Lambda>0$ such that
$\Rs+h\in{\cal R}$ whenever $\,\|h\|_{L^2(Q)}\le\lambda\,$ for some $\lambda\in (0,\Lambda]$.
In the following, we shall only consider such small perturbations $h$. Let, for any such $h\in L^2(Q)$,
\begin{align*}
&\vph={\cal S}(\Rs+h),\quad \muh=-\Delta\vph+f'(\vph), \quad \uh=-\nabla\ph+\muh\nabla\vph,
\end{align*}
where ${\rm mean}\,\ph(t)=0$, $\aat$, and
\begin{align*}
&\yh=\vph-\vps-\xi^h, \quad \zh=\muh-\mus-\eta^h,\quad \wh=\uh-\bus-\bv^h,
\quad \rh=\ph-\ps-q^h.
\end{align*}
Observe that $\,{\rm mean}\,\rh(t)=0$ for a.e. $t\in (0,T)$ as well.

According to the definition of the notion of Fr\'echet differentiability,
it suffices to show the existence of an increasing function $\,Z:(0,\Lambda)
\mapsto (0,+\infty)\,$ such that $$\lim_{\lambda\searrow0}\,\frac{Z(\lambda)}{\lambda^2}=0$$
 and
\Beq\label{Frechet}
\|\yh\|_{{\cal W}}^2\,\le\,Z(\|h\|_{L^2(Q)}) \quad\mbox{for all }\,h\in L^2(Q) \quad\mbox{with }
\,\|h\|_{L^2(Q)}\le\Lambda\,.
\Eeq

We are going to show that we may choose $\,Z(\lambda)=\widetilde C\,\lambda^4\,$ with some suitably
chosen constant $\widetilde C>0$.
To this end, we first observe that, according to Theorem \ref{exe2D} and Lemma \ref{line}, we have,
for all admissible variations $\,h\,$
the following regularity properties for $(\yh, \zh, \wh, \rh)$:
\begin{align}
\label{reghvar}
 \yh&\in H^1(0,T;\Ldue)\cap C^0([0,T];\Hdue)\cap L^2(0,T;H^4(\Omega)),\non\\[1mm]
\zh&\in C^0([0,T];\Ldue)\cap L^2(0,T;\Hdue),\non\\[1mm]
\wh&\in L^2(0,T;L^4(\Omega)),\quad \rh\in L^2(0,T;W^{2,4/3}(\Omega)).
\end{align}
Moreover, we see from Theorem \ref{exe2D} and Lemma \ref{conti} that there hold the bounds \eqref{ssbounds1} and \eqref{ssbounds2}
for both $(\vps,\mus,\bus,p^*)$ and $(\vph,\muh,\uh,\ph)$, as well as the stability estimate
\begin{align}\label{stabu3}
&\|\vph-\vps\|_{H^1(0,t;\Ldue)\cap C^0([0,t];\Hdue)\cap L^2(0,t;H^4(\Omega))}\,+\,\|\mu^h-\mus\|_{L^2(0,t;\Hdue)}\non\\
&\qquad +\|\uh-\bus\|
_{L^2(0,t;L^4(\Omega))}\,+\,\|\ph-p^*\|_{L^2(0,t;W^{2,4/3}(\Omega))}\non\\
&\quad \le\,K_2\,\|h\|_{L^2(0,t;\Ldue)}\, ,
\end{align}
for some $K_2>0$ depending only on $\|\vp_0\|_{\Hdue}$, $\Omega$, $T$, $\|S\|_{L^2(Q)}$ and $\widehat{R}$.
Besides, simple
algebraic manipulations, together with the facts $${\rm div}(\uh)={\rm div}(\bus)=S,\quad {\rm div}({\bf v}^h)=0,$$
yield that, by its definition,  $(\yh,\zh,\wh,\rh)$ is a strong solution to the following problem on $Q$:
\begin{align}
\label{h1}
&\dt\yh-\Delta \zh\non\\
&\quad =\,-\yh\,S-\wh\cdot\nabla\vps
-\bus\cdot \nabla \yh-(\uh-\bus) \cdot\nabla(\vph-\vps) \quad\aeQ,\\[1mm]
\label{h2}
&\zh\,=\,-\Delta \yh +f'(\vph)-f'(\vps)-f''(\vps)\xi^h \quad\aeQ,\\[1mm]
\label{h3}
&\wh\,=\,-\nabla\rh + \zh\nabla\vps+\mus\nabla\yh+(\muh-\mus)\nabla(\vph-\vps)
\quad\aeQ,\\[1mm]
\label{h4}
&{\rm div}(\wh)=0\quad\aeQ,\\[1mm]
\label{h5}
&-\Delta\rh\,=\,-{\rm div}(\zh\nabla\vps)-{\rm div}(\mus\nabla\yh)-{\rm div}((\muh-\mus)\nabla(\vph-\vps))\quad\aeQ,\\[1mm]
\label{h6}
&\dn\yh=\dn\zh=\dn \rh=0 \quad\mbox{and }\,\wh\cdot{\bf n}=0, \quad\aeS,\\[1mm]
\label{h7}
&\yh(0)=0 \quad\aeO.
\end{align}

We now perform a number of estimates for the system \eqref{h1}--\eqref{h7}, where
in the remainder of this proof the letters $\,C\,$ and $\,C_i$, $i\in\enne$, stand for positive constants
that may depend on the data but not on the choice of $h\in L^2(Q)$ that satisfies
$\Rs+h\in{\cal R}$ and $\|h\|_{L^2(Q)}\,\le\,\Lambda$. The actual value
of $\,C\,$ may change within lines or even within formulas.
To begin with, we recall Taylor's formula
$$f'(\vph)=f'(\vps)+f''(\vps)(\vph-\vps)+\frac12 f^{(3)}(\sigma^h)(\vph-\vps)^2,$$
where $\,\sigma^h=a\vph+(1-a)\vps\,$ with some $\,a(x,t)\in [0,1]\,$ for $\,(x,t)\in Q$. Then,
 by the definition of $y^h$, we have that
\Beq
\label{Taylor1}
\left|f'(\vph)-f'(\vps)-f''(\vps)\xi^h\right|
\,\le\,C\left|\yh\right|\,+\,C\left|\vph-\vps\right|^2
\quad\aeQ.
\Eeq
We observe that \eqref{h2} yields for almost every $t\in (0,T)$ that
$$
\left|{\rm mean}\,\zh(t)\right|\,\le\,C\iO\left(\left|\yh(t)\right|\,+\,\left|\vph(t)-\vps(t)\right|^2\right)\,dx,
$$
whence
\begin{align*}
\int_0^t\left|{\rm mean}\, \zh(s)\right|^2\,ds\,\le\,C\iQt\left|\yh\right|^2\,dxds\,+\,C\int_0^t\left\|\vph(s)-\vps(s)\right\|^4_{L^2(\Omega)}\,ds\,.
\end{align*}
Thus, it follows from Poincar\'e's inequality and \eqref{stabu3} that
\begin{align}
\label{Poin2}
&\int_0^t\left\|\zh(s)\right\|^2_{\Huno}\,ds\non\\
&\quad \le\,C\iQt\left|\nabla\zh\right|^2\,dxds\,+\,C\int_0^t\left\|\yh(s)\right\|^2_{\Ldue}\,ds\,+\,
C\,\|h\|_{L^2(Q_t)}^4\,.
\end{align}

\vspace{2mm}\noindent
{\bf First estimate:}

\vspace{1mm}\noindent
We multiply \eqref{h1} by $\,\yh$, account for \eqref{h2}, and integrate over $Q_t$, where $t\in (0,T]$. After
integration by parts, we then obtain the identity
\begin{align}
\label{fd1}
&\frac 12\,\left\|\yh(t)\right\|_{\Ldue}^2\,+\iQt\left|\Delta\yh\right|^2\,dxds\non\\
&\quad =\,\iQt \left(f'(\yh)-f'(\vps)-f''(\vps)\xi^h\right)\Delta\yh\, dxds\non\\
&\qquad -\iQt\yh\left(\yh\,S+\wh\cdot\nabla\vps+\bus\cdot \nabla \yh\right)dxds\non\\
&\qquad -\iQt\yh\,\nabla(\vph-\vps)\cdot(\uh-\bus)\,dxds\non\\
&\quad :=\,\sum_{j=1}^3I_j\,,
\end{align}
with obvious notation.  From \eqref{Taylor1} and \eqref{stabu3}, we obtain for every $\gamma\in(0,1)$ (to be specified later) that
\begin{align}
\label{fd2}
|I_1|\,&\le\,\gamma\iQt\left|\Delta\yh\right|^2\,dxds\,+\,\frac C\gamma\iQt\left|\yh\right|^2\,dxds+\,\frac C\gamma
\iQt\left|\vph-\vps\right|^4\,dxds\non\\[1mm]
&\le\,\gamma\iQt\left|\Delta\yh\right|^2\,+\,\frac C\gamma\iQt\left|\yh\right|^2\,+\,\frac C\gamma\,\|h\|_{L^2(Q_t)}^4\,.
\end{align}
Next, by virtue of $\mathrm{div}(\bus)=S$, H\"older's inequality and Young's inequality, we also have
\begin{align}
\label{fd3}
|I_2|\,&=  \left|\iQt\yh\left(\frac12 \yh\,S+\wh\cdot\nabla\vps\right)\,dxds\right|\non\\
&\le\int_0^t\left\|\yh(s)\right\|_{L^4(\Omega)}\left(\left\|\yh(s)\right\|_{L^4(\Omega)}\,\|S(s)\|_{\Ldue}\,+\,\left\|\wh(s)\right\|_{\Ldue}\,
\|\nabla\vps(s)\|_{L^4(\Omega)}\right)ds\non\\
&\le\,\frac 18\iQt\left|\wh\right|^2\,dxds\,+\,C\int_0^t(1+\|S(s)\|_{\Ldue}^2)\left\|\yh(s)\right\|_{\Huno}^2\,ds\,.
\end{align}
Moreover, using \eqref{stabu3} once more, we see that
\begin{align}
\label{fd4}
|I_3|\,&\le\int_0^t\left\|\yh(s)\right\|_{\Ldue}\,\left\|\nabla(\vph(s)-\vps(s))\right\|_{L^4(\Omega)}\,\left\|\uh(s)-\bus(s)\right\|_{L^4(\Omega)}\,ds
\non\\[1mm]
&\le\,\left\|\vph-\vps\right\|_{C^0([0,t];\Hdue)}\,\left\|\uh-\bus\right\|_{L^2(0,t;L^4(\Omega))}\,
\left\| \yh \right\|_{L^2(Q_t)}\non\\
&\le \iQt\left|\yh\right|^2\,dxds\,+\,C\,\|h\|_{L^2(Q_t)}^4\,.
\end{align}
Combining \eqref{fd1}--\eqref{fd4}, we have thus shown the estimate
\begin{align}
\label{fd5}
& \left\|\yh(t)\right\|_{\Ldue}^2\,+\,(2-2\gamma)\iQt\left|\Delta\yh\right|^2\,dxds\non\\
&\quad \le\,\frac 14\iQt\left|\wh\right|^2 dxds\,+\,C\left(1+\gamma^{-1}\right)\|h\|_{L^2(Q_t)}^4\non\\
&\qquad +\,C\left(1+\gamma^{-1}\right)\int_0^t(1+\|S(s)\|_{\Ldue}^2)\left\|\yh(s)\right\|^2_{\Huno}\,ds\,.
\end{align}

\vspace{2mm}\noindent
{\bf Second estimate:}

\vspace{1mm}\noindent
We now multiply \eqref{h1} by $\zh$, take the scalar product of \eqref{h3} with $\wh$, add the
resulting identities, and integrate over $Q_t$. After integration by parts, we then get
\begin{align}
\label{fd6}
&\frac 12\,\left\|\nabla\yh(t)\right\|_{\Ldue}^2\,+\iQt\!\left|\nabla\zh\right|^2\,dxds\,+\,\iQt\!\left|\wh\right|^2\,dxds\non\\
&\quad =\,
-\!\iQt\! \yh\,S\,\zh\,dxds\,-\,\iQt \zh \bus\cdot\nabla \yh\,dxds\,-\,\iQt\!\zh\,(\uh-\bus)\cdot \nabla(\vph-\vps)\,dxds\non\\
&\qquad +\,\iQt\!\mus\,\nabla\yh \cdot \wh\,dxds\,+\,\iQt\!(\muh-\mus)\nabla(\vph-\vps)\cdot \wh\,dxds\non\\
&\qquad -\iQt\dt\yh\left[f'(\vph)-f'(\vps)-f''(\vps)\xi^h\right]\,dxds \non\\
&\quad :=\,\,\sum_{j=1}^6J_j\,,
\end{align}
with obvious notation. We have, by the H\"older and Young inequalities and \eqref{Poin2}, that
\begin{align}
\label{fd7}
|J_1|\,&\le\int_0^t\|S(s)\|_{\Ldue}\left\|\zh(s)\right\|_{L^4(\Omega)}\left\|\yh(s)\right\|_{L^4(\Omega)} \, ds\non\\
&\le\,\gamma\iQt\left|\nabla\zh
\right|^2\,dxds\,+\,C\gamma \|h\|_{L^2(Q_t)}^4\non\\
&\quad +\,C\left(\gamma+\gamma^{-1}\right)\int_0^t\!\left(1+\|S(s)\|_{\Ldue}^2\right)\left\|\yh(s)\right\|^2_{\Huno}\,ds\,.
\end{align}
Moreover, by similar reasoning, we obtain
\begin{align}
\label{fd7a}
|J_2| & \le \int_0^t \|\zh(s)\|_{L^4(\Omega)}\|\bus(s)\|_{L^4(\Omega)}\|\nabla \yh\|_{\Ldue}\, ds\non\\
&\le\,\gamma\iQt\left|\nabla\zh
\right|^2\,dxds\,+\,C\gamma\|h\|_{L^2(Q_t)}^4\non\\
&\quad +\,C\left(\gamma+\gamma^{-1}\right)\int_0^t\!\left(1+\|\bus(s)\|_{L^4(\Omega)}^2\right)\left\|\yh(s)\right\|^2_{\Huno}\,ds\,,
\end{align}
and
\begin{align}
\label{fd8}
|J_3|\,&\le\int_0^t\left\|\zh(s)\right\|_{\Ldue} \left\|\uh(s)-\bus(s)\right\|_{L^4(\Omega)} \left\|\nabla(\vph(s)-\vps(s))\right\|_{L^4(\Omega)}\,ds
\non\\
&\le\,C\left\|\vph-\vps\right\|_{C^0([0,t];\Hdue)}\left\|\zh\right\|_{L^2(Q_t)}\left\|\uh-\bus\right\|
_{L^2(0,t;L^4(\Omega))}\non\\
&\le\,\gamma\iQt\left|\nabla\zh\right|^2\,dxds\,+\, C(\gamma+\gamma^{-1})\,\|h\|_{L^2(Q_t)}^4\,+\,C\gamma \int_0^t\left\|\yh(s)\right\|^2_{\Ldue}\,ds\,.
\end{align}
We also have, invoking \eqref{Lady} and standard elliptic estimates, that
\begin{align}
\label{fd9}
|J_4|\,&\le\int_0^t\|\mus(s)\|_{L^4(\Omega)} \left\|\nabla\yh(s)\right\|_{L^4(\Omega)} \left\|\wh(s)\right\|_{\Ldue}\,ds\non\\
&\le\,\frac 14\iQt\left|\wh\right|^2\,dxds\non\\
&\quad +\,C\int_0^t \|\mus(s)\|_{L^2(\Omega)} \|\mus(s)\|_{H^1(\Omega)}\left\|\nabla\yh(s)\right\|_{\Ldue}\left\|\nabla\yh(s)\right\|_{\Huno}\,ds\non\\
&\le\,\frac 14\iQt\left|\wh\right|^2\,dxds\,+\,\gamma\iQt\left|\Delta\yh\right|^2\,dxds\non\\
&\quad +\,C(1+\gamma^{-1})\int_0^t \left(1+\|\mus(s)\|_{H^1(\Omega)}^2\right)\left\|\yh(s)\right\|^2_{\Huno}\,ds\,.
\end{align}
Besides,
\begin{align}\label{fd10}
|J_5|\,&\le\int_0^t \left\|\muh(s)-\mus(s)\right\|_{L^4(\Omega)} \left\|\nabla(\vph(s)-\vps(s))\right\|_{L^4(\Omega)} \left\|\wh(s)\right\|_{\Ldue} \,ds
\non\\
&\le\,C\,\left\|\vph-\vps\right\|_{C^0([0,t];\Hdue)}\left\|\muh-\mus\right\|
_{L^2(0,t;\Huno)}\,\left\|\wh\right\|_{L^2(Q_t)} \non\\
&\le\,\frac 14\iQt\left|\wh\right|^2\,dxds\,+\,C\,\|h\|_{L^2(Q_t)}^4\,.
\end{align}
The estimation of $J_6$ requires some preparations. First notice that almost everywhere in $Q$ it holds
\begin{align*}
&f'(\vph)-f'(\vps)-f''(\vps)\xi^h\non\\
&\quad =\,\int_0^1\frac d{ds}\left(f'(s\,\vph+(1-s)\,\vps)\right)ds\,-\,f''(\vps)\xi^h\\
&\quad =\,\left(\vph-\vps\right)\int_0^1 f''(s\,\vph+(1-s)\,\vps)\,ds\,-\,f''(\vps)\xi^h\\
&\quad =\,f''(\vps)\,\yh\,+\left(\vph-\vps\right)\int_0^1\bigl[ f''(s\,\vph+(1-s)\,\vps)-f''(\vps)\bigr]\,ds\\
&\quad =\,f''(\vps)\,\yh\,+\left(\vph-\vps\right)^2\int_0^1\!\!\int_0^1 s\,f^{(3)}\left(\sigma(s\,\vph+(1-s)\,\vps)+(1-\sigma)\,\vps
\right) d\sigma\,ds\\[1mm]
&\quad :=\,f''(\vps)\,\yh\,+\,A^h\,,
\end{align*}
with obvious notation. Thus, using the estimate \eqref{ssbounds2} for $\vph$ and $\vps$, we have almost everywhere on $Q$ the estimates
\begin{align}
\label{fd11}
\left|A^h\right|\,&\le\,C\left|\vph-\vps\right|^2, \\[1.5mm]
\label{fd12}
\left|\dt A^h\right|\,&\le\,C\left|\vph-\vps\right|\left|\dt\vph-\dt\vps\right|\,+\,C\left|\vph-\vps\right|^2
\left(\left|\dt\vph\right|+\left|\dt\vps\right|\right)\,.
\end{align}
Now, we write $J_6$ as
\Beq\label{fd13}
J_6\,=\,-\iQt f''(\vps)\,\yh\,\dt\yh\,dxds\,-\,\iQt\dt\yh\,A^h dxds\,:=\,J_7+J_8,
\Eeq
with obvious notation.
Clearly,
$$J_7\,=\,-\frac 12 \iQt\dt\left(f''(\vps)\left|\yh\right|^2\right)\,dxds\,+\,\frac 12\iQt\left|\yh\right|^2\,f^{(3)}(\vps)\,
\dt\vps\,dxds,$$
that is, since $\,\,-f''(\vps)=1-3{\vps}^2\,\le\,1$,
\begin{align}\label{fd14}
|J_7|\,&\le\,\frac 12\left\|\yh(t)\right\|^2_{\Ldue}\,+\,C\int_0^t\left\|\dt\vps(s)\right\|_{\Ldue}\left\|\yh(s)\right\|_{L^4(\Omega)}^2\,ds
\non\\
&\le\,\frac 12\left\|\yh(t)\right\|^2_{\Ldue}\,+\,C\int_0^t\left\|\dt\vps(s)\right\|_{\Ldue}\left\|\yh(s)\right\|_{\Huno}^2\,ds\,.
\end{align}
Moreover,
$$J_8\,=\,-\iQt\dt\!\left(\yh\,A^h\right)\,dxds\,+\,\iQt \yh\,\dt A^h\,dxds,$$
so that, owing to \eqref{fd12},
\begin{align}
\label{fd15}
|J_8|\,&\le\iO\left|\yh(t)\right|\left|A^h(t)\right|\,dx\,+\,C\iQt\left|\yh\right| \left|\vph-\vps\right|
\left|\dt\vp^h-\dt\vp^*\right|\,dxds \non\\
&\quad +\, C\iQt \left|\yh\right| \left|\vph-\vps\right|^2\left(\left|\dt \vp^h\right|+\left|\dt \vp^*\right|\right)\,dxds\non\\
&\le\,\frac 14 \left\| \yh(t) \right\|_{\Ldue}^2\,+\,\,C\,\left\|\vph(t)-\vps(t)\right\|_{L^4(\Omega)}^4\non\\
&\quad\,\,+\,C\int_0^t\left\|\yh(s)\right\|_{L^4(\Omega)}\left\|\vph(s)-\vps(s)\right\|_{L^4(\Omega)}
\left\|\dt \vp^h(s)- \dt\vp^*(s)\right\|_{\Ldue}\,ds\non\\
&\quad\,\, +\,C\int_0^t\left\|\yh(s)\right\|_{L^6(\Omega)}\left\|\vph(s)-\vps(s)\right\|_{L^6(\Omega)}^2
\left(\left\|\dt \vp^h(s)\right\|_{\Ldue}\,+\,\left\|\dt \vp^*(s)\right\|_{\Ldue}\right)ds\non\\
&\le\,\frac 14\left\|\yh(t)\right\|_{\Ldue}^2\,+\,C\,\|h\|_{L^2(Q_t)}^4\non\\[1mm]
&\quad\,\,+\,C\left\|\vph-\vps\right\|_{C^0([0,t];\Huno)}\left\|\dt \vp^h-\dt \vp^*\right\|_{L^2(Q_t)}\left\|
\yh\right\|_{L^2(0,t;\Huno)}\non\\[1mm]
&\quad\,\,+\,C\left\|\vph-\vps\right\|_{C^0([0,t];\Huno)}^2\left(\left\|\dt \vp^h\right\|_{L^2(Q_t)}
+\left\|\dt \vp^*\right\|_{L^2(Q_t)}\right)\left\|
\yh\right\|_{L^2(0,t;\Huno)}\non\\
&\le\,\frac 14\left\|\yh(t)\right\|_{\Ldue}^2\,+\,C\,\|h\|_{L^2(Q_t)}^4\,+\,C\int_0^t\left\|\yh(s)\right\|^2_{\Huno}\,ds\,.
\end{align}
Combining the estimates \eqref{fd6}--\eqref{fd15}, we obtain the following inequality:
\begin{align}
\label{fd16}
&\frac 12\left\|\nabla\yh(t)\right\|^2_{\Ldue}\,+\,(1-3\gamma)\iQt\left|\nabla\zh\right|^2\,dxds\,
+\,\frac 12\iQt\left|\wh\right|^2\,dxds\non\\
&\quad \leq \frac 34\left\|\yh(t)\right\|_{L^2(\Omega)}^2\,+\,\gamma \iQt |\Delta \yh|^2\, dxds\,+\,C(1+\gamma+\gamma^{-1})\,\|h\|_{L^2(Q_t)}^4\non\\
&\qquad + C\left(1+\gamma+\gamma^{-1}\right)\int_0^t\Psi_4(s)\left\|\yh(s)\right\|^2_{\Huno}\,ds
\,,
\end{align}
where the function
$$\Psi_4(s)= 1+\|S(s)\|^2_{\Ldue}+\|\bus(s)\|_{\Ldue}^2+\|\mus(s)\|_{\Huno}^2+\left\|\dt\vp^*(s)\right\|_{\Ldue}^2$$ is known to belong to $L^1(0,T)$. Hence,
adding the inequalities \eqref{fd5} and \eqref{fd16},  we obtain
\begin{align}
\label{fd16a}
& \frac14\left(\left\|\nabla\yh(t)\right\|^2_{\Ldue}+\left\|\yh(t)\right\|_{\Ldue}^2\right)\,+\,(2-3\gamma)\iQt\left|\Delta\yh\right|^2\,dxds\non\\
&\qquad +\,(1-3\gamma)\iQt\left|\nabla\zh\right|^2\,dxds\,
+\,\frac 14\iQt\left|\wh\right|^2\,dxds\non\\
&\quad \le C\left(1+\gamma+\gamma^{-1}\right)\|h\|_{L^2(Q_t)}^4\,+\,C\left(1+\gamma+\gamma^{-1}\right)\int_0^t \Psi_4(s)\left\|\yh(s)\right\|^2_{\Huno}\,ds\,.
\end{align}
Adjusting $\gamma\in(0,1)$ small enough and applying Gronwall's lemma,
using also the estimate \eqref{Poin2}, we deduce that
\begin{align*}
\left\|\yh\right\|^2_{C^0([0,t];\Huno)\cap L^2(0,t;\Hdue)}\,+\,\left\|\zh\right\|^2_{L^2(0,t;\Huno)}
\,+\,\left\|\wh\right\|^2_{L^2(Q_t)}\,\le\,C_2\,\|h\|^4_{L^2(Q_t)}\,,
\end{align*}
for all $t\in(0,T]$. The condition \eqref{Frechet} is therefore fulfilled with the choice $\,Z(\lambda)=C_2\lambda^4$. As a consequence, we see that
$$\frac{\|\mathcal{S}(R^*+h)-\mathcal{S}(R^*)-\xi^h\|_{\mathcal{W}}}{\|h\|_{L^2(Q)}}=\frac{\|y^h\|_{\mathcal{W}}}{\|h\|_{L^2(Q)}}
\leq \sqrt{C_2}\,\|h\|_{L^2(Q)}\to 0$$
as $\|h\|_{L^2(Q)}\to 0$. This concludes the proof of the assertion.
\Edim

%%%%%%%%%%%%%%%%%%%%%%%%%%%%%%%%%%%%%%%%%%%%%%%%%%%%%%%%%%%%%%%%%%%%%%%%

\section{The Optimal Control Problem}
\setcounter{equation}{0}
In this section, we establish our main results for the optimal control problem {\bf (CP)} stated in the introduction.
In addition to the previous assumptions (A1)--(A3), we assume the following

\vspace{1.5mm}\noindent
(A4) \quad The target functions satisfy $\vp_\Omega\in L^2(\Omega)$, $\vp_Q\in L^2(Q)$. The constants $\beta_i$, $i=1,2,3$,\linebreak \hspace*{13mm}
are nonnegative but not  all zero, and the functions $R_{\rm min}, R_{\rm max}
\in L^\infty(Q)$ satisfy \linebreak
\hspace*{13mm} $R_{\rm min}\le R_{\rm max}$ almost everywhere in $Q$.

\subsection{Existence}

We begin with the result on the existence of an optimal control.
\Bthm\label{opexe}
Suppose that the assumptions {\rm (A1)--(A4)} are satisfied. Then the optimal control problem ${\bf (CP)}$ admits at least one solution $(\vp, R)$ such that
$R\in \Uad$ and $\vp=\mathcal{S}(R)$ is the
unique global strong solution to problem \eqref{ssphi}--\eqref{ssini}.
\Ethm
\Bdim
The proof essentially follows from the convexity/coercivity of the nonnegative
cost functional $\mathcal{J}(\vp, R)$. Consider the reduced cost functional
\begin{align}
\widetilde{\mathcal{J}}(R): L^2(Q)\to [0,+\infty) \quad\text{such that}\quad \widetilde{\mathcal{J}}(R):=\mathcal{J}(\vp, R),\label{ReJ}
\end{align}
for any $R\in L^2(Q)$, where $\vp=\mathcal{S}(R)$ is the corresponding
unique global strong solution to problem \eqref{ssphi}--\eqref{ssini} with given initial datum $\vp_0$ and source $S$ satisfying (A1), (A2).
We observe that problem ${\bf (CP)}$ is equivalent to the minimization problem
$$\min_{R\in \Uad} \widetilde{\mathcal{J}}(R).$$
Then we pick a bounded minimizing sequence $\{R_n\}_{n\in\enne}$, i.e., a sequence of admissible controls
such that $\lim_{n\to+\infty}\widetilde{\mathcal{J}}(R_n)=\inf_{R\in \Uad} \widetilde{\mathcal{J}}(R)$ and put
\Beq
\label{ex1}
\vp_n\,=\,{\cal S}(R_n),\quad \mu_n\,=\,-\Delta\vp_n+f'(\vp_n),\quad \bu_n\,=\,-\nabla p_n+\mu_n\nabla\vp_n,
\Eeq
where $\,p_n\,$ is the unique solution to the elliptic boundary value problem
resembling \eqref{ssp} with the constraint $\,{\rm mean}\,p_n(t)=0\,$ for almost every $t\in (0,T)$.
Since $\Uad$ is a bounded closed subset of $L^2(Q)$ (see (A4)), in view
of Theorem \ref{exe2D}, in particular \eqref{ssbounds1}, we may infer from standard compactness
arguments (cf. \cite[Sect.~8, Cor.~4]{Simon}) that there are some $R\in\Uad$ and a quadruple $(\vp,\mu,\bu,p)$
such that, at least for some subsequence which is again indexed by $n$ for simplicity,
\begin{align*}
R_n\to R\quad&\mbox{weakly-star in $\,L^\infty(Q)$},\\[1mm]
\vp_n\to\vp\quad&\mbox{weakly-star in $H^1(0,T;\Ldue)\cap L^\infty(0,T;\Hdue)\cap L^2(0,T;H^4(\Omega))$}
\\
&\mbox{strongly in $C^0([0,T];W^{1,r}(\Omega))$ for $1\le r<+\infty$, in $C^0(\overline Q)$},\\[1mm]
&\mbox{and also in $L^2(0,T; H^3(\Omega))$},\\[1mm]
\mu_n\to\mu\quad&\mbox{weakly-star in $L^\infty(0,T;\Ldue)\cap L^2(0,T;\Hdue)$},\\[1mm]
\bu_n\to\bu\quad&\mbox{weakly in $L^2(0,T;\Huno)^2$},\\[1mm]
p_n\to p\quad&\mbox{weakly in $L^2(0,T;\Hdue)$}.
\end{align*}
It follows, in particular, that $\vp,\mu,\bu$ satisfy \eqref{ssdivu}--\eqref{ssini}. Moreover, we can deduce the following convergence results for nonlinear terms:
\begin{align*}
f'(\vp_n)\to f'(\vp)\quad &\mbox{strongly in $C^0(\overline Q)$},\\[1mm]
\mu_n\nabla\vp_n \to\mu\nabla\vp\quad & \mbox{weakly in $L^2(Q)^2$},\\[1mm]
{\rm div}(\vp_n\bu_n)\to{\rm div}(\vp\bu)\quad& \mbox{weakly in $L^2(Q)$},\\[1mm]
{\rm div}(\mu_n\nabla\vp_n)\to{\rm div}(\mu\nabla\vp)\quad& \mbox{weakly in $L^2(0,T; (\Huno)')$}.
\end{align*}
These convergent results enable us to pass to the limit as $n\to\infty$ in the state system \eqref{ssphi}--\eqref{ssini}, \eqref{ssp} (at least in the weak formulation), and
noting that the constraint \eqref{mean} also holds true. Thus, we infer that $(\vp,\mu,\bu,p)$ is the (unique) solution to the state system associated with the control $R$, namely, it holds  $\vp=\mathcal{S}(R)$.
Therefore, the limit $(\vp, R)$ is an admissible pair for the control problem {\bf (CP)}. It then follows from the sequential
semicontinuity properties of the cost functional $\mathcal{J}$ that $(\vp,R)$ is a solution to the control problem {\bf (CP)} with $R$ being an optimal control.
\Edim

\subsection{First-order necessary optimality conditions}

With the Fr\'echet differentiability for $\mathcal{S}$ that has been shown in Theorem \ref{diffS}, we can conclude from the convexity
of $\Uad$ and the chain rule of differentiation the following first-order necessary optimality condition for problem {\bf (CP)}.
\Bthm\label{noc}
Suppose that the assumptions {\rm (A1)--(A4)} are fulfilled. Let $R^*\in\Uad$ be a solution to the optimal control
problem {\bf (CP)} with the associated state $\vps={\cal S}(R^*)$. Then we have the following variational inequality
\begin{align}\label{vug1}
& \beta_1\iO\left(\vps(T)-\vp_\Omega\right)\xi(T)\,dx\,+\,\beta_2\int_Q\left(\vps-\vp_Q\right)\xi\,dxdt\non\\
& \qquad +\,\beta_3 \int_Q R^*(R-R^*)\,dxdt\,\ge\,0
\end{align}
for all $R\in\Uad$, where $\,\xi\,$ is the first component of the unique solution to the linearized system
\eqref{lsphi}--\eqref{ellq} with $\,h=R-R^*$.
\Ethm
\Bdim
Recalling the definition of the reduced cost functional $\widetilde{\mathcal{J}}$ (see \eqref{ReJ}) and invoking the convexity of $\Uad$, we obtain (cf. \cite[Lemma 2.21]{To})
$$\big(\widetilde{\mathcal{J}}'(R^*), R-R^*\big)\geq 0\,\quad \forall\, R\in \Uad.$$
By the chain rule, we have $\widetilde{\mathcal{J}}'(R)= \mathcal{J}'_{\mathcal{S}(R)}(\mathcal{S}(R), R)\circ D\mathcal{S}(R)+\mathcal{J}'_{R}(\mathcal{S}(R), R)$, where for every fixed $R\in L^2(Q)$,
$\mathcal{J}'_{\vp}(\vp, R)$ is the Fr\'echet derivative of $\mathcal{J}(\vp,R)$  with respect to $\vp$ at $\vp\in \mathcal{W}$ and for every fixed $\vp\in \mathcal{W}$, $\mathcal{J}'_{R}(\vp, R)$ is the Fr\'echet derivative with respect to $R$ at $R\in L^2(Q)$. Then, by a straightforward computation and using the fact $D\mathcal{S}(R^*)(R-R^*)=\xi$ (see Theorem \ref{diffS}), we obtain the variational inequality \eqref{vug1}.
\Edim

We now turn our interest to the derivation of first-order necessary optimality conditions, where
we aim at expressing the two summands in the variational inequality \eqref{vug1} containing the component $\,\xi\,$ (i.e., the solution to the linearized system \eqref{lsphi}--\eqref{ellq})
in terms of the adjoint state variables.

To this end, assume that the assumptions (A1)--(A4) are satisfied and that
a fixed optimal control $R^*\in\Uad$ is given. Let $(\vps,\mus,\bus,p^*)$ denote the corresponding unique solution to the state system
\eqref{ssphi}--\eqref{ssini}, \eqref{ssp} established in Theorem \ref{exe2D}. Moreover, let $(\xi,\eta,\bv,q)$ be the
unique solution to the linearized system \eqref{lsphi}--\eqref{ellq}, according to Lemma \ref{line}.
We then consider the following \emph{adjoint system}: find a quadruple $(p_1, p_2,{\bf p}_3, p_4)$
such that
\begin{align}
\label{asphi}
&-\dt p_1-\bus\cdot\nabla p_1+\Delta p_2-f''(\vps)\,p_2+\mathbf{p}_3 \cdot \nabla \mus\,=\,\beta_2(\vps-\vp_Q)\quad\mbox{a.e. in \,$Q$,}\\[1mm]
\label{asmu}
& p_2\,=\,\Delta p_1\,+\, \mathbf{p}_3\cdot \nabla \vps\quad\mbox{a.e. in \,$Q$},\\[1mm]
\label{asu}
& {\bf p}_3\,=\nabla p_4\,-\,p_1\nabla\vps\quad\mbox{a.e. in \,$Q$,}\\[1mm]
\label{asp}
&\mathrm{div}(\mathbf{p}_3)=0\quad\mbox{a.e. in \,$Q$,}
\end{align}
with the following boundary and endpoint conditions:
\begin{align}
\label{asbc}
&\dn p_1\,=\,\dn p_2\,=\,{\bf p}_3\cdot{\bf n}\,=\,0\quad\mbox{a.e. on \,$\Sigma$,}\\[1mm]
\label{asini}
&p_1(T)\,=\,\beta_1\,(\vps(T)-\vp_\Omega)\quad\mbox{a.e. in \,$\Omega$.}
\end{align}
Besides, we see from \eqref{asu}--\eqref{asbc}  that $p_4$ satisfies
\begin{align}
&\Delta p_4\,=\,{\rm div}(p_1\nabla\vps)\quad\mbox{a.e. in \,$Q$,}\label{asp4}\\[1mm]
&\dn p_4=0\quad \mbox{a.e. on $\Sigma$}.\label{asp4bc}
\end{align}
Observe that $p_4$ is only determined up to a constant. Thus, we make $p_4$ unique by postulating that
$\,{\rm mean}\,p_4(t)=0\,$ for almost every $t\in (0,T)$.

The adjoint system \eqref{asphi}--\eqref{asini} turns out to be a backward-in-time problem and it can be easily derived by using
the formal Lagrange method described in \cite{To} with direct computations via integration by parts.
Owing to the fact that we only have $p_1(T)=\beta_1(\vps(T)-\vp_\Omega)\in \Ldue$ (recall (A4)), the system \eqref{asphi}--\eqref{asini} in general cannot
be expected to enjoy a strong solution on $Q$ (unless $\beta_1=0$, or $\vps(T)=\vp_\Omega$), and the corresponding solution
$(p_1, p_2, \mathbf{p}_3, p_4)$ can only be expected to have a certain weaker regularity.
Therefore, instead of the pointwise equations \eqref{asphi}--\eqref{asmu}, $(p_1, p_2, \mathbf{p}_3, p_4)$ should be understood as a solution satisfying the weak formulation:
\begin{align}
\label{wp1}
&\langle -\dt p_1,\psi\rangle_{(\Hdue)',\Hdue}\,+\iO\Delta p_1\,\Delta\psi
\,dx + \iO (\mathbf{p}_3\cdot \nabla \vps)\Delta \psi\,dx\non\\[1mm]
&\qquad - \iO (f''(\vps)\,\Delta p_1+\bus\cdot\nabla p_1)\,\psi\,dx\,-\, \iO f''(\vps)(\mathbf{p}_3\cdot \nabla \vps)\psi\,dx \non\\
&\qquad + \iO (\mathbf{p}_3\cdot \nabla \mus) \psi\, dx\non\\
&\quad =\,\iO\beta_2(\vps-\vp_Q)\,\psi\,dx,\qquad\mbox{for all $\,\psi\in H^2(\Omega)\,$ and a.e. in $(0,T)$}.%
\end{align}

We have the following result on the existence of adjoint states.
\Blem\label{adexe}
Suppose that the assumptions {\rm (A1)--(A4)} are fulfilled and that the optimal control $R^*$ as well as the associate state $(\vps,\mus,\bus,p^*)$ are given as above. Then the adjoint problem
\eqref{asu}--\eqref{wp1} admits a unique weak solution $(p_1, p_2, \mathbf{p}_3, p_4)$ on $Q$ satisfying
\begin{align}
& p_1 \in H^1(0,T;(\Hdue)')\cap C^0([0,T];\Ldue)\cap L^2(0,T;\Hdue), \label{wreg1}\\[1mm]
& p_2\in L^2(Q), \quad{\bf p}_3\in L^2(0,T; H^1(\Omega)^2), \quad p_4\in L^2(0,T;\Hdue).\label{wreg2}
\end{align}
\Elem
\Bdim
The proof follows from a similar argument as that used in the proof of Lemma \ref{line}, namely, by
means of the Faedo--Galerkin procedure. Therefore, we simply omit the implementation of the approximation
scheme and just perform the necessary \textit{a priori} estimates. Taking $\psi=p_1$ in the weak form \eqref{wp1}, we have
\begin{align}
&-\frac12\frac{d}{dt}\|p_1\|_{\Ldue}^2 \,+\,\|\Delta p_1\|_{\Ldue}^2\non\\
&\quad = - \iO (\mathbf{p}_3\cdot \nabla \vps)\Delta p_1\,dx\,+\,\iO \big[f''(\vps)\,\Delta p_1+\bus\cdot\nabla p_1\big]\,p_1\,dx\non\\[1mm]
&\qquad + \iO f''(\vps)(\mathbf{p}_3\cdot \nabla \vps)p_1 \,dx \, -\, \iO (\mathbf{p}_3\cdot \nabla \mus)p_1 \, dx\non\\[1mm]
&\qquad +\,\iO\beta_2(\vps-\vp_Q)\,p_1\,dx\non\\
&\quad :=\sum_{j=1}^5 I_j.
\label{esline}
\end{align}
Recalling the uniform estimate \eqref{ssbounds1} for $\vps$, we have
\begin{align}
\|\mathrm{div}(p_1\nabla \vps)\|_{L^{4/3}(\Omega)}
&\leq \,C\,\|\nabla p_1\|_{\Ldue} \, \|\nabla \vps\|_{L^4(\Omega)}\,+\, C\,\|p_1\|_{L^4(\Omega)}\,\|\Delta \vps\|_{\Ldue}\non\\
&\leq C\,\| p_1\|_{\Huno}.\non
\end{align}
Then similar to \eqref{stime3}, it holds that
\begin{align}
\|\nabla p_4\|_{L^4(\Omega)}
&\leq \,C\,\|p_4\|_{W^{2,4/3}(\Omega)}\,\leq\, C\,\|\mathrm{div}(p_1\nabla \vps)\|_{L^{4/3}(\Omega)}\non\\
&\leq  \,C\,\| p_1\|_{\Huno}.\label{p4L4}
\end{align}
On the other hand, it follows from H\"older's inequality that
\begin{align}
\|p_1\nabla\vps\|_{L^4(\Omega)}\leq C\,\|p_1\|_{L^8(\Omega)}\,\|\nabla\vps\|_{L^8(\Omega)}\leq C\,\| p_1\|_{\Huno},\non
\end{align}
which together with \eqref{p4L4} yields
\begin{align}
\|\mathbf{p}_3\|_{L^4(\Omega)}\leq C\,\| p_1\|_{\Huno}.\label{p3L4}
\end{align}
As a consequence, it follows that
\begin{align}
|I_1|&\leq \|\mathbf{p}_3\|_{L^4(\Omega)}\,\|\nabla \vps\|_{L^4(\Omega)}\,\|\Delta p_1\|_{\Ldue}\non\\
&\leq C\,\| p_1\|_{\Huno}\,\|\Delta p_1\|_{\Ldue}\non\\
&\leq \frac18\,\|\Delta p_1\|_{\Ldue}^2 + C\,\| p_1\|_{\Ldue}^2.\label{adI1}
\end{align}
Next, using the fact $\mathrm{div}(\bus)=S$, we have
\begin{align}
|I_2|&= \left| \iO \left(f''(\vps)\,(\Delta p_1)p_1 \,+\, \frac12 \bus\cdot  \nabla ( p_1)^2 \right)\,dx  \right|\non\\
&= \left| \iO \left(f''(\vps)\,(\Delta p_1)p_1 \,-\, \frac12 S (p_1)^2\right) \,dx  \right|\non\\
&\leq \|f''(\vps)\|_{L^\infty(\Omega)} \, \|\Delta p_1\|_{\Ldue} \, \|p_1\|_{\Ldue}+\|S\|_{\Ldue}\,\|p_1\|_{L^4(\Omega)}^2\non\\
&\leq \frac18 \,\|\Delta p_1\|_{\Ldue}^2 + C\,(1+\|S\|_{\Ldue}^2)\,\| p_1\|_{\Ldue}^2.\label{adI2}
\end{align}
Besides, we have
\begin{align}
|I_3|&\leq \|f''(\vps)\|_{L^\infty(\Omega)}\,\|\mathbf{p}_3\|_{L^4(\Omega)}\,\|\nabla \vps\|_{L^4(\Omega)}\,\|p_1\|_{\Ldue}\non\\[2mm]
&\leq C\,\|p_1\|_{\Huno}\,\|p_1\|_{\Ldue}\non\\
&\leq \frac18\,\|\Delta p_1\|_{\Ldue}^2 + C\,\| p_1\|_{\Ldue}^2,\label{adI3}
\end{align}
\begin{align}
|I_4|&\leq \|\mathbf{p}_3\|_{L^4(\Omega)}\,\|\nabla \mus\|_{\Ldue}\,\|p_1\|_{L^4(\Omega)}\non\\[2mm]
&\leq C\,\|\nabla \mus\|_{\Ldue}\,\|p_1\|_{\Huno}^2\non\\
&\leq \frac18\,\|\Delta p_1\|_{\Ldue}^2 + C\,(1+\|\nabla \mus\|_{\Ldue}^2)\,\| p_1\|_{\Ldue}^2,\label{adI4}
\end{align}
and
\begin{align}
|I_5|& \leq \beta_2\|\vps-\vp_Q\|_{\Ldue}\,\|p_1\|_{\Ldue}\leq \|p_1\|_{\Ldue}^2+  \frac{\beta_2^2}{4}\, \|\vps-\vp_Q\|_{\Ldue}^2.\label{adI5}
\end{align}
From \eqref{esline} and \eqref{adI1}--\eqref{adI5}, we arrive at the following differential inequality
\begin{align}
&-\frac{d}{dt}\|p_1\|_{\Ldue}^2 \,+\,\|\Delta p_1\|_{\Ldue}^2 \leq C\,\Psi_5(t)\,\| p_1\|_{\Ldue}^2 + \beta_2^2\, \|\vps-\vp_Q\|_{\Ldue}^2
\end{align}
for a.e. $t\in (0,T)$, where $\Psi_5(t)=1+\|S\|_{\Ldue}^2+\|\nabla \mus\|_{\Ldue}^2\in L^1(0,T)$. Then, by the (backward) Gronwall inequality, we obtain
\begin{align}
&\|p_1(t)\|_{\Ldue}^2+ \int_t^T\|\Delta p_1(s)\|_{\Ldue}^2\,ds\non\\
&\quad  \leq C\Big(\|p_1(T)\|_{\Ldue}^2 + \beta_2^2 \|\vps-\vp_Q\|_{L^2(Q)}^2\Big),\quad \forall\, t\in [0,T].
\end{align}

 The above estimate yields that $p_1\in L^\infty(0,T; L^2(\Omega))\cap L^2(0,T; H^2(\Omega))$. Recalling \eqref{p4L4} and \eqref{p3L4},
 we also infer that $p_4\in L^2(0,T; W^{2,4/3}(\Omega))$ and $\mathbf{p}_3\in L^2(0,T; L^4(\Omega))$. Then,
 by a comparison argument in \eqref{wp1}, we obtain $\dt p_1\in L^2(0,T; (H^2(\Omega))')$, which further implies that\, $p_1\in C^0([0,T]; \Ldue)$.
Since
\begin{align}
\|\mathrm{div}(p_1\nabla \vps)\|_{\Ldue}&\leq C\,\|\nabla p_1\|_{L^4(\Omega)}\,\|\nabla \vps\|_{L^4(\Omega)}\,+\,
C\,\|p_1\|_{L^\infty(\Omega)}\,\|\Delta \vps\|_{L^2(\Omega)}\non\\
&\leq C\,\|p_1\|_{H^2(\Omega)}\in L^2(0,T),\non
\end{align}
then, by the standard elliptic estimate, we infer that $p_4\in L^2(0,T; H^2(\Omega))$. This fact and \eqref{asu}
yield that\, $\mathbf{p}_3\in L^2(0,T; H^1(\Omega)^2)$.

Thus, we are able to prove the existence of a weak solution to the adjoint problem
 \eqref{asu}--\eqref{wp1} satisfying the regularity properties \eqref{wreg1}--\eqref{wreg2}. Besides, for this linear system, the proof of uniqueness is straightforward, and we omit the details here.
\Edim

Now we are able to eliminate the function $\xi$ from the variational inequality
\eqref{vug1} and, alternatively, form a first-order necessary optimality condition by the
state system \eqref{ssphi}--\eqref{ssini} for $\vps$ together with the adjoint system \eqref{asu}--\eqref{wp1}:

\Bcor\label{noc1}
Suppose that the assumptions {\rm (A1)--(A4)}  are fulfilled, and let $R^*\in\Uad$ be a solution to the optimal control
problem {\bf (CP)} with the associated state $\vps={\cal S}(R^*)$ as well as
the adjoint state $p_1$. Then we have the following variational inequality:
\begin{align}
\label{vug2}
\int_Q\left(p_1\,+\,\beta_3\,R^*\right)\left(R-R^*\right)\,dxds \ge\,0\, \qquad\forall\,R\in\Uad.
\end{align}
\Bdim
We infer from \eqref{lsini}, \eqref{asini} and the Newton--Leibniz formula that
\begin{align}
&\beta_1\iO\left(\vps(T)-\vp_\Omega\right)\xi(T)\,dx\, =\, \int_0^T \frac{d}{dt}\left(\int_\Omega p_1 \xi\,dx\right) dt.\label{be1}
\end{align}
On the other hand, taking $(\xi, \eta, \mathbf{v}, q)$, which is the unique solution to the linearized system
\eqref{lsphi}--\eqref{ellq} with $\,h=R-R^*$, as test functions in the adjoint system \eqref{asu}--\eqref{wp1}, adding the results together and using integration by parts,
we have
\begin{align}
&\beta_2\int_Q\,\left(\vps-\vp_Q\right)\xi\,dxdt\non\\
&\quad = -\int_0^T\,\langle\dt p_1, \xi\rangle_{(\Hdue)',\Hdue} dt\,+\,\int_Q p_2\, \Delta \xi\,dxdt\non\\
&\qquad -\int_Q \Big[\bus\cdot\nabla p_1-f''(\vps)\,p_2+\mathbf{p}_3 \cdot \nabla \mus\Big]\,\xi\,dxdt\non\\
&\qquad +\int_Q ( p_2\,-\,\Delta p_1\,-\, \mathbf{p}_3\cdot \nabla \vps)\,\eta\,dxdt\non\\
&\qquad +\int_Q ({\bf p}_3\,-\,\nabla p_4\,+\,p_1\nabla\vps)\cdot \mathbf{v}\,dxdt\,-\,\int_Q q\, \mathrm{div}(\mathbf{p}_3)\, dxdt\non\\
&\quad = - \int_0^T \frac{d}{dt}\left(\int_\Omega p_1 \xi\,dx\right) dt\, +\, \int_Q \Big[\dt \xi\, + \,{\rm div}(\vps\,\bv)\,+\,{\rm div}(\xi\,\bus)\,-\,\Delta\eta\Big]\, p_1\, dxdt\non\\
&\qquad  + \int_Q \Big[\eta\,+\,\Delta\xi\,-\,f''(\vps)\xi\Big]\, p_2\,dxdt\, +\, \int_Q p_4\, \mathrm{div}(\mathbf{v})\,dxdt\non\\
&\qquad  + \int_Q(\bv \,+\,\nabla q\,-\,\eta\,\nabla\vps\,-\,\mus\nabla\xi) \cdot \mathbf{p}_3\,dxdt\non\\
&\quad = - \int_0^T \frac{d}{dt}\left(\int_\Omega p_1 \xi\,dx\right) dt\,+\, \int_Q (R-R^*)\, p_1\,dxdt,
\label{be2}
\end{align}
where in the last step, we have used the equations \eqref{lsphi}--\eqref{lsdivu} for $(\xi, \eta, \mathbf{v}, q)$. Adding \eqref{be1} and \eqref{be2} together, we obtain
$$\beta_1\iO\left(\vps(T)-\vp_\Omega\right)\xi(T)\,dx\, +\, \beta_2\int_Q\left(\vps-\vp_Q\right)\xi\,dxdt \,=\, \int_Q (R-R^*) p_1\,dxdt.$$
Hence, the variational inequality \eqref{vug2} is an immediate consequence of the above identity and \eqref{vug1}.
\Edim
\Brem
When $\beta_3>0$, it follows from \eqref{vug2} that the optimal control $\,R^*\,$ is nothing but the $L^2(Q)$-orthogonal projection of
\,$-\beta_3^{-1}p_1\,\,$ onto the closed convex set $\Uad$. Then, by a standard argument, we infer the following pointwise condition:
\begin{align}
R^*(x,t)=\max\left\{R_{\mathrm{min}}(x,t),\ \min\{-\beta_3^{-1}p_1(x,t),\,R_{\mathrm{max}}(x,t)\}\right\},\quad \text{for a.e.}\ (x,t)\in Q.\non
\end{align}
\Erem
\Ecor
%%%%%%%%%%%%%%%%%%%%%%%%%%%%%%%%%%%%%%%%%%%%%%%%%%%%%%%%%%%%%%%%%%%%%%%%

\section*{Acknowledgements}
 The research of H. Wu is partially supported by NNSFC grant No. 11631011 and the Shanghai Center for Mathematical Sciences.

%\section{Complements}

%\setcounter{equation}{0}

%%%%%%%%%%%%%%%%%%%%%%%%%%%%%%%%%
%% bibliography
%%%%%%%%%%%%%%%%%%%%%%%%%%%%%%%%%

\vspace{3truemm}

%%%%%%%%%%%%%%%%%%%%%%%%%%%%%%%%%
%% bibliography
%%%%%%%%%%%%%%%%%%%%%%%%%%%%%%%%%

\vspace{3truemm}

\Begin{thebibliography}{10}
\itemsep=0pt
\bibitem{Anderson}
D.-M. Anderson, G.-B. McFadden and A.-A. Wheeler,
Diffuse-interface methods in fluid mechanics,
{\it Annual Review of Fluid Mech.} {\bf 30}(1) (1997), 139--165.

\bibitem{BLM}
N. Bellomo, N.-K. Li and P.-K. Maini,
On the foundations of cancer modelling: Selected topics, speculations, and perspectives,
\textit{Math. Models Methods Appl. Sci.} \textbf{18} (2008), 593--646.

\bibitem{BCG}
S. Bosia, M. Conti and M. Grasselli,
On the Cahn--Hilliard--Brinkman system,
\textit{Commun. Math. Sci.} \textbf{13} (2015), 1541--1567.

\bibitem{Cahn}
J.\,W. Cahn and J.-E. Hilliard,
Free energy of a nonuniform system. I. Interfacial free energy,
{\it J. Chem. Phys.} {\bf 28} (1958), 258--267.

\bibitem{CWSL}
Y. Chen, S.-M. Wise, V.-B. Shenoy, J.-S. Lowengrub,
A stable scheme for a nonlinear multiphase tumor growth model with an elastic membrane,
\textit{Int. J. Numer. Methods Biomed. Eng.} \textbf{30} (2014), 726--754.

\bibitem{CFGS1}
P. Colli, M.-H. Farshbaf-Shaker, G. Gilardi and J. Sprekels,
Optimal boundary control of a viscous Cahn--Hilliard system with dynamic boundary condition and double obstacle potentials,
\textit{SIAM J. Control Optim.} \textbf{53} (2015), 2696--2721.

\bibitem{CGS3}
P. Colli, G. Gilardi and J. Sprekels,
A boundary control problem for the viscous Cahn--Hilliard equation with dynamic boundary conditions,
\textit{Appl. Math. Optim.} \textbf{72} (2016), 195--225.

\bibitem{CGS4}
P. Colli, G. Gilardi and J. Sprekels,
Optimal velocity control of a viscous Cahn--Hilliard system with convection and dynamic boundary conditions,
\textit{SIAM J. Control Optim.} \textbf{56} (2018), 1665--1691.

\bibitem{CGS5}
P. Colli, G. Gilardi and J. Sprekels,
Optimal velocity control of a convective Cahn--Hilliard system with double obstacles and dynamic boundary conditions: a `deep quench' approach, to appear in \textit{J. Convex Anal.} \textbf{26} (2019) (see also Preprint arXiv:1709.03892 [math. AP] (2017), 1--30).

\bibitem{CGH15}
P. Colli, G. Gilardi and D. Hilhorst,
On a Cahn--Hilliard type phase field system related to tumor growth,
\textit{Discrete Contin. Dyn. Syst.} \textbf{35} (2015), 2423--2442.

\bibitem{CGRS}
P. Colli, G. Gilardi, E. Rocca and J. Sprekels,
Vanishing viscosities and error estimate for a Cahn--Hilliard type phase field system related to tumor growth,
\textit{Nonlinear Anal. Real World Appl.} \textbf{26} (2015), 93--108.

\bibitem{CGRS1}
P. Colli, G. Gilardi, E. Rocca and J. Sprekels,
Optimal distributed control of a diffuse interface model of tumor growth,
\textit{Nonlinearity} \textbf{30} (2017), 2518--2546.

\bibitem{CGS1}
P. Colli, G. Gilardi and J. Sprekels,
A boundary control problem for the pure Cahn--Hilliard equation with dynamic boundary conditions,
\textit{Adv. Nonlinear Anal.} \textbf{4} (2015), 311--325.

\bibitem{CG}
M. Conti and A. Giorgini,
The three-dimensional Cahn--Hilliard--Brinkman system with unmatched viscosities,
preprint, 2018.
https://hal.archives-ouvertes.fr/hal-01559179

\bibitem{CL2010}
V. Cristini and J.-S. Lowengrub,
\textit{Multiscale Modeling of Cancer: An Integrated Experimental and Mathematical Modeling Approach},
Cambridge Univ. Press, Cambridge, 2010.

\bibitem{DFRSS}
M. Dai, E. Feireisl, E. Rocca, G. Schimperna and M. Schonbek,
Analysis of a diffuse interface model for multi-species tumor growth,
\textit{Nonlinearity} \textbf{30}(4) (2017), 1639.

\bibitem{DG}
F. Della Porta and M. Grasselli,
On the nonlocal Cahn--Hilliard--Brinkman and Cahn--Hilliard--Hele--Shaw systems,
\textit{Commun. Pure Appl. Anal.} \textbf{15} (2016), 299--317,
Erratum: \textit{Commun. Pure Appl. Anal.} \textbf{16} (2017), 369--372.

\bibitem{DGG}
F. Della Porta, A. Giorgini and M. Grasselli,
The nonlocal Cahn--Hilliard--Hele--Shaw system with logarithmic potential,
\textit{Nonlinearity} \text{31} (2018), 4851.

\bibitem{FBG2006}
A. Fasano, A. Bertuzzi and A. Gandolfi,
Mathematical modeling of tumour growth and treatment,
Complex Systems in Biomedicine, Springer, Milan, 71--108, 2006.

\bibitem{FW2012}
X. Feng and S.-M. Wise,
Analysis of a Darcy--Cahn--Hilliard diffuse interface model for the Hele--Shaw flow and its fully discrete finite element approximation,
\textit{SIAM J. Numer. Anal.} \textbf{50} (2012), 1320--1343.

\bibitem{Fri2007}
A. Friedman,
Mathematical analysis and challenges arising from models of tumor growth,
\textit{Math. Models Methods Appl. Sci.} \textbf{17} (2007), 1751--1772.

\bibitem{Lowen10}
H.-B. Frieboes, F. Jin, Y.-L. Chuang, S.-M. Wise, J.-S. Lowengrub and V. Cristini,
Three-dimensional multispecies nonlinear tumor growth - II: tumor invasion and angiogenesis,
{\it J. Theoret. Biol.} {\bf 264} (2010), 1254--1278.

\bibitem{FGR}
S. Frigeri, M. Grasselli and E. Rocca,
On a diffuse interface model of tumor growth,
\textit{European J. Appl. Math.} \textbf{26} (2015), 215--243.

\bibitem{FGS}
S. Frigeri, M. Grasselli and J. Sprekels,
Optimal distributed control of two-dimensional nonlocal Cahn--Hilliard--Navier--Stokes systems with degenerate mobility and singular potential, \textit{Appl. Math. Optim.}, online first article, 2018. DOI: 10.1007/s00245-018-9524-7.

\bibitem{FLRS}
S. Frigeri, K.-F. Lam, E. Rocca and G. Schimperna,
On a multi-species Cahn--Hilliard--Darcy tumor growth model with singular potentials,
\textit{Comm Math Sci.} \textbf{16} (2018), 821--856.

\bibitem{FRS}
S. Frigeri, E. Rocca and J. Sprekels,
Optimal distributed control of a nonlocal Cahn--Hilliard/Navier--Stokes system in two dimensions,
\textit{SIAM J. Control. Optim.} \textbf{54} (2016), 221--250.

\bibitem{GL2016}
H. Garcke and K.-F. Lam,
Global weak solutions and asymptotic limits of a Cahn--Hilliard--Darcy system modelling tumour growth,
\textit{AIMS Mathematics}, \textbf{1}(3) (2016), 318--360.

\bibitem{GL2018}
H. Garcke and K.-F. Lam,
On a Cahn--Hilliard--Darcy system for tumour growth with solution dependent source terms,
In: Rocca E., Stefanelli U., Truskinovsky L., Visintin A. (eds),
Trends in Applications of Mathematics to Mechanics, Springer INdAM Series, Vol \textbf{27}, Springer, 2018.

\bibitem{GLNS}
H. Garcke, K.-F. Lam, R. N\"urnberg and E. Sitka,
A multiphase Cahn--Hilliard--Darcy model for tumour growth with necrosis,
\textit{Math. Models Methods Appl. Sci.} \textbf{28}(3) (2018), 525--577.

\bibitem{GLR}
H. Garcke, K.-F. Lam and E. Rocca,
Optimal control of treatment time in a diffuse interface model for tumour growth,
\textit{Appl. Math. Optim.}, online first article, 2017. DOI: 10.1007/s00245-017-9414-4.

\bibitem{GLSS}
H. Garcke, K.-F. Lam, E. Sitka and V. Styles,
A Cahn--Hilliard--Darcy model for tumour growth with chemotaxis and active transport,
\textit{Math. Models Methods Appl. Sci.} \textbf{26}(6) (2016), 1095--1148.

\bibitem{GioGrWu}
A. Giorgini, M. Grasselli and H. Wu,
The Cahn--Hilliard--Hele--Shaw system with singular potential,
{\it Ann. Inst. H. Poincare Anal. Non Lineaire}, {\bf 35}(4) (2018), 1079--1118.

\bibitem{HZO12}
A. Hawkins-Daarud, K.-G. van der Zee and J.-T. Oden,
Numerical simulation of a thermodynamically consistent four-species tumor growth model,
\textit{Int. J. Numer. Meth. Biomed. Engrg.} \textbf{28} (2012), 3--24.

\bibitem{HKW17}
H. Hinterm\"uller, M. Keil and D. Wegner,
Optimal control of a semi-discrete Cahn--Hilliard--Navier--Stokes system with non-matched fluid densities,
{\it SIAM J. Control Optim.} {\bf 55} (2017), 1954--1989.

\bibitem{HiWe12}
M. Hinterm\"uller and D. Wegner,
Distributed optimal control of the Cahn--Hilliard system including the case of a double-obstacle homogeneous free energy density,
\textit{SIAM J. Control Optim.} \textbf{50} (2012), 388--418.

\bibitem{HiWe14}
M. Hinterm\"uller and D. Wegner,
Optimal control of a semi-discrete Cahn--Hilliard--Navier--Stokes system,
\textit{SIAM J. Control Optim.} \textbf{52} (2014), 747--772.

\bibitem{HiWe17}
M. Hinterm\"uller and D. Wegner,
Distributed and boundary control problems for the semi-discrete Cahn--Hilliard/Navier--Stokes system with non-smooth Ginzburg--Landau energies,
Topological Optimization and Optimal Transport, Radon Series on Computational and Applied Mathematics, \textbf{17} (2017), 40--63.

\bibitem{JWZ}
J. Jiang, H. Wu and S. Zheng,
Well-posedness and long-time behavior of a non-autonomous
Cahn--Hilliard--Darcy system with mass source modeling tumor growth,
{\it J. Differential Equations} {\bf 259} (2015), 3032-3077.

\bibitem{LTZ}
J.-S. Lowengrub, E.-S. Titi and K. Zhao,
Analysis of a mixture model of tumor growth,
{\it European J. Appl. Math.} \textbf{24} (2013), 691--734.

\bibitem{LT}
J. Lowengrub and L. Truskinovsky,
Quasi-incompressible Cahn--Hilliard fluids and topological transitions,
{\it R. Soc. Lond. Proc. Ser. A Math. Phys. Eng. Sci.} {\bf 454}(1978) (1998), 2617--2654.

\bibitem{RS}
E. Rocca and J. Sprekels,
Optimal distributed control of a nonlocal convective Cahn--Hilliard equation by the velocity in three dimensions,
\textit{SIAM J. Control Optim.} \textbf{53} (2015), 1654--1680.

\bibitem{Simon}
J. Simon,
Compact sets in the space $L^p(0,T; B)$,
{\it Ann. Mat. Pura Appl.~(4)\/} {\bf 146} (1987), 65--96.

\bibitem{To}
F. Tr\"oltzsch,
Optimal Control of Partial Differential Equations. Theory, Methods and Applications,
Grad. Stud. in Math., Vol. \textbf{112}. AMS, Providence, RI, 2010.

\bibitem{WW2012}
X.-M. Wang and H. Wu,
Long-time behavior for the Hele--Shaw--Cahn--Hilliard system,
\textit{Asymptot. Anal.} \textbf{78} (2012), 217--245.

\bibitem{WZ2013}
X.-M. Wang and Z.-F. Zhang,
Well-posedness of the Hele--Shaw--Cahn--Hilliard system,
\textit{Ann. Inst. H. Poincar\'e Anal. Non Lin\'eaire}, \textbf{30} (2013), 367--384.

\bibitem{Wise2010}
S.-M. Wise,
Unconditionally stable finite difference, nonlinear multigrid simulations of the Cahn--Hilliard--Hele--Shaw system of equations,
\textit{J. Sci. Comput.} \textbf{44} (2010), 38--68.

\bibitem{Wise2011}
S.-M. Wise, J.-S. Lowengrub and V. Cristini,
An adaptive multigrid algorithm for simulating solid tumor growth using mixture models,
\textit{Math. Comput. Modelling} \textbf{53} (2011), 1--20.

\bibitem{Lowen08}
S.-M. Wise, J.-S. Lowengrub, H.-B. Frieboes and V. Cristini,
Three dimensional multispecies nonlinear tumor growth - I: model and numerical method,
{\it J. Theoret. Biol.} {\bf 253} (2008), 524--543.

\bibitem{ZL1}
X.-P. Zhao and C.-C. Liu,
Optimal control of the convective Cahn--Hilliard equation,
\textit{Appl. Anal.} \textbf{92} (2013), 1028--1045.

\bibitem{ZL2}
X.-P. Zhao and C.-C. Liu,
Optimal control of the convective Cahn--Hilliard equation in 2D case,
\textit{Appl. Math. Optim.} \textbf{70} (2014), 61--82.

\End{thebibliography}

\End{document}

%%%%%%%%%%%%%%%%%%%%%%%%%%%%%%%%%%%%%%%%%%%%%